\documentclass[10pt]{article}

\usepackage{amssymb}

\newtheorem{theorem}{Theorem}
\newtheorem{definition}[theorem]{Definition}
\newtheorem{lemma}[theorem]{Lemma}
\newtheorem{proposition}[theorem]{Proposition}
\newtheorem{corollary}[theorem]{Corollary}

\newcommand{\var}{\mbox{\rm var}}																									
\newcommand{\C}{\cal C}
\newcommand{\B}{\cal B}
\newcommand{\A}{\cal A}

\def\qed{\unskip\nobreak\hfill\penalty50\hskip 3pt\hbox{}\nobreak
\hfill\hbox{\vrule width 4 pt height 10 pt}}

\begin{document}

\title{The Central Limit Theorem for uniformly strong mixing measures}
\author{Nicolai Haydn\thanks{Mathematics Department, USC,
Los Angeles, 90089-1113. E-mail: $<$nhaydn@math.usc.edu$>$.
This work was supported by a grant from the NSF (DMS-0301910).}}
\maketitle

\begin{abstract}
The theorem of Shannon-McMillan-Breiman states that for every generating
partition on an ergodic system,
the exponential decay rate of the measure of cylinder sets 
equals the metric entropy almost everywhere (provided the entropy is finite).
In this paper we prove that the measure of cylinder sets are lognormally 
distributed for strongly mixing systems and infinite partitions and show that the rate of convergence
is polynomial provided the fourth moment of the information function is finite. 
Also, unlike previous results by Ibragimov and others which only apply to finite partitions,
here we do not require any regularity of the conditional entropy function.
We also obtain the law of the iterated logarithm and the weak invariance principle for the 
information function. 
\end{abstract}

\section{Introduction}

Let $\mu$ be a $T$-invariant probability measure on a space $\Omega$ on 
which the map $T$ acts measurably. For a  measurable partition $\cal A$
one forms the $n$th join ${\cal A}^n=\bigvee_{j=0}^{n-1}T^{-j}{\cal A}$
which forms a finer partition of $\Omega$. (The atoms of ${\cal A}^n$ are
traditionally called {\em $n$-cylinders}.)  For $x\in\Omega$
we denote by $A_n(x)\in{\cal A}^n$ the $n$-cylinder which contains $x$.
The Theorem of Shannon-McMillan-Breiman (see e.g.~\cite{Man,Petersen}) then states that for 
$\mu$-almost every $x$ in $\Omega$ the limit
$$
\lim_{n\rightarrow\infty}\frac{-\log\mu(A_n(x))}n
$$
exists and equals the metric entropy $h(\mu)$ provided the entropy is finite in the case
of a countable infinite partition. It is easy to see that this
convergence is not uniform (not even for Bernoulli measures with 
weights that are not all equal). This theorem was proved for finite partitions in increasing degrees
of generality in the years 1948 to 1957 and then was by Carleson~\cite{Car}
 and Chung~\cite{Chu}
generalised to infinite partitions. Similar results (for finite partitions) for the recurrence
and waiting times were later proved by Ornstein and Weiss~\cite{OW} and
Nobel and Wyner~\cite{NW} respectively. The limiting behaviour for recurrence
times was generalised in 2002 by Ornstein and Weiss~\cite{OW2} to countably infinite
partitions.  In the present paper we are concerned
with the limiting distribution of the information function $I_n(x)=-\log\mu(A_n(x))$
around its mean value.

The statistical  properties of $I_n$ are of great interest in information theory where they are
connected to the efficiency of compression schemes. Let us also note that in dynamical systems
the analog of SMB's theorem for compact metric spaces is the Brin-Katok local entropy formula
\cite{BK}
which states that for  ergodic invariant measures the exponential decay rate of dynamical 
balls is almost everywhere equal to the entropy. 

There is a large classical body of work on the Central Limit Theorem (CLT)
for independent random variables. For dependent random 
variables the first CLTs are due to Markov (for Markov chains) and 
Bernstein~\cite{Ber} for 
random variables that are allowed to have some short range dependency but
have to be independent if separated by a suitable time difference (for more
than a power of the length $n$ of the partial sums $S_n$). 
In 1956 Rosenblatt~\cite{R1} then introduced the notions of uniform mixing and 
strong mixing (see below) and proved a  CLT for the
partial sums $S_n$ of random variables that satisfy the strong mixing property. 
In \cite{R2} he then proved a more general CLT for random variables on systems
that satisfy an $L^2$ norm condition\footnote{
The map $T$ satisfies an {\em $L^2$ norm condition}  if 
$\displaystyle\sup_{f:\mu(f)=0}\frac{\|T^nf\|^2}{\|f\|^2}$
decays exponentially fast as $n\rightarrow\infty$. This is a somewhat stronger
mixing condition than the strong mixing condition
}.
Around the same time Nagaev~\cite{Nag} proved a convergence
theorem for the stable law for strongly mixing systems. His result
covers the case of the CLT and formed the basis for Ibragimov's famous
1962 paper~\cite{Ibr} in which he proved for finite partitions
 `a refinement to SMB's theorem' by 
 showing that  $I_n(x)=-\log\mu(A_n(x))$
is in the limit lognormally distributed for systems that are strongly 
mixing and satisfy a regularity condition akin to a Gibbs property. 
Based on his results and methods, Philipp and Stout~\cite{PS}
proved the almost sure invariance principle for the information function
 $I_n$ under similar conditions as Ibragimov used (requiring faster
decay rates). This result in turn was then used by
Kontoyiannis~\cite{Kon} to prove the almost sure invariance principle,
CLT and the law of the iterated logarithm LIL for recurrence and waiting times, 
thus strengthening the result
of Nobel and Wyner~\cite{NW} who showed that for strongly mixing systems (without
regularity condition) the exponential growth rate of waiting times equals the
metric entropy.

Various improvements and refinements to the CLT for the 
information function have been
successively done mainly for measures that satisfy a genuine Gibbs 
property. For instance Collet, Galves and Schmitt~\cite{CGS} in 
order to prove the lognormal distribution of entry times for 
exponentially $\psi$-mixing Gibbs measures\footnote{
We say an invariant probability measures $\mu$ is {\em Gibbs}
for a potential $f$ with pressure $P(f)$
if there exists a constant $c>0$ so that 
$$
\frac1c\le\frac{\mu(A_n(x))}{e^{f(x)+f(Tx)+\cdots+f(T^{n-1}x)-nP(f)}}\le c
$$
for every $x\in \Omega$ and $n=1,2,\dots$.
}
 needed to know
that $I_n$ is in the limit lognormally distributed. A more general 
result is due to Paccaut~\cite{Pac} for maps on the interval where
he had to assume some topological covering properties. For some
non-uniformly hyperbolic maps on the interval similar results
were formulated in  \cite{FHV, BV}. However all those results use
explicitly the Gibbs property of the invariant measure $\mu$ 
 to approximate the information function $I_n$ by an ergodic sum
and then to invoke standard results on the CLT for 
sufficiently regular observables (see for instance~\cite{Gor,Liv,ch}). 
 (Of course the variance has to be 
non-zero because otherwise the limiting distribution might not be normal 
as an example in \cite{CGS} illustrates.) 

Results
that do not require the explicit Gibbs characterisation of the measure 
like Kontoyiannis' paper~\cite{Kon}, all ultimately rely on the original 
paper of Ibragimov~\cite{Ibr} and require apart from the strong mixing 
condition the regularity of the Radon-Nikodym derivative of the measure
under the local inverse maps. 
In \cite{HV2} we went beyond his regularity 
constraint and proved a CLT with error bounds for the lognormal distribution of
the information function for $(\psi,f)$-mixing systems which 
included traditional $\psi$-mixing maps and also equilibrium states
for rational maps with critical points in the Julia set.

The present paper is significant in two respects: (i) we allow 
for the partition to be countably infinite instead of finite and (ii) unlike Ibragimov
(and all who followed him) we do not require
an $L^1$-regularity condition for the Radon Nikodym derivative for
local inverses of the map. This condition which was introduced in \cite{Ibr}
 is the $L^1$ equivalent of
what otherwise would allow a transfer operator approach to analyse
the invariant measure and imply the Gibbs property\footnote{
More precisely, Ibragimov's condition requires that the $L^1$-norms of the
differences $f-f_n$ decay polynomially, where $f=\lim_{n\rightarrow\infty}f_n$
and $f_n(x)=\log\mathbb{P}(x_0|x_{-1}x_{-2}\dots x_{-n})$.
}. 
We moreover prove that the rate of convergence is polynomial (Theorem~\ref{CLT})
and the variance is always positive for genuinely infinite partitions.

Let us note that convergence rates for the CLT have previously been
obtained by A Broise~\cite{Bro} for a large class of expanding maps on the interval for
which the Perron-Frobenius operator has a `spectral gap'. Similar 
estimates were obtained by P\`ene~\cite{Pen} for Gibbs measures for
dispersing billiards.

This paper is structured as follows:
In the second section we introduce uniform strong mixing systems and in the
third section we
prove the existence of the variance $\sigma^2$ of strongly mixing probability 
measures (Proposition~\ref{variance.convergence}) as well as 
 the growth rate of higher order moments (Proposition~\ref{fourth.moment}). 
 This is the main part of the proof (note that Ibragimov's regularity condition
 was previously needed precisely to obtain the variance of the measure).
In section~4 we then prove the CLT using Stein's method of exchangeable
pairs. In section~5 we prove the Weak Invariance Principle for $I_n$ using the 
CLT and the convergence rate obtained in section~3. 

I would like to thank my colleague Larry Goldstein for many conversations in which 
he explained Stein's method to me.

\section{Main results}
Let $T$ be a map on a space $\Omega$ and $\mu$ a
probability measure on $\Omega$. Moreover let $\cal A$ be a (possibly infinite)
 measurable partition of $\Omega$ and denote by
${\cal A}^n=\bigvee_{j=0}^{n-1}T^{-j}{\cal A}$
its {\em $n$-th join} which also is a measurable partition of $\Omega$ for
every $n\geq1$. The atoms of ${\cal A}^n$ are called {\em $n$-cylinders}.
Let us put ${\cal A}^*=\bigcup_{n=1}^\infty{\cal A}^n$ for the collection of
all cylinders in $\Omega$ and put $|A|$ for the length of a
cylinder $A\in{\cal A}^*$, i.e.\ $|A|=n$ if $A\in{\cal A}^n$.

We shall assume that $\cal A$ is generating, i.e.\ that the atoms of
${\cal A}^\infty$ are single points in $\Omega$.

\subsection{Mixing}

\begin{definition}\label{definition}
We say the invariant probability measure $\mu$ is {\em uniformly strong mixing}
if there exists a decreasing function $\psi: \mathbb{N}\rightarrow\mathbb{R}^+$ 
which satisfies $\psi(\Delta)\rightarrow0$ as $\Delta\rightarrow\infty$ so that
$$
\left|\sum_{(B,C)\in S}\left(\mu(B\cap C)-\mu(B)\mu(C)\right)\right|\le\psi(\Delta)
$$
for every subset $S$ of ${\cal A}^n\times T^{-\Delta-n}{\cal A}^m$ and every $n,m,\Delta>0$.
\end{definition}

\noindent {\bf Various kinds of mixing:}\footnote{Here we adopt probabilistic terminology 
which differs from the one used in the dynamical systems community.}\\
In the following list of different mixing properties $U$ is always in the $\sigma$-algebra
 generated by ${\cal A}^n$ and $V$ lies in the $\sigma$-algebra 
generated by ${\cal A}^*$ (see also~\cite{Dou}).
The limiting behaviour is as the length of the `gap' $\Delta\rightarrow\infty$:
\begin{enumerate}
\item  {\em $\psi$-mixing}: 
$\displaystyle
\sup_n\sup_{U,V}\left|\frac{\mu(U\cap T^{-\Delta-n}V)}{\mu(U)\mu(V)}-1\right|\rightarrow0.
$
\item  {\em Left $\phi$-mixing}:
$\displaystyle
\sup_n\sup_{U,V}\left|\frac{\mu(U\cap T^{-\Delta-n}V)}{\mu(U)}-\mu(V)\right|\rightarrow0.
$
\item {\em Strong mixing} \cite{R1,Ibr} (also called $\alpha$-mixing):
$\displaystyle
\sup_n\sup_{U,V}\left|\mu(U\cap T^{-\Delta-n}V)-\mu(U)\mu(V)\right|\rightarrow0$.
\item   {\em Uniform mixing}~\cite{R1,R2}:
$\displaystyle
\sup_n\sup_{U,V}\left|\frac1k\sum_{j=1}^k\mu(U\cap T^{-n-j}V)-\mu(U)\mu(V)\right|\rightarrow0
$ as $k\rightarrow\infty$.

\end{enumerate}
One can also have {\em right $\phi$-mixing} 
when
$
\sup_n\sup_{U,V}\left|\frac{\mu(U\cap T^{-\Delta-n}V)}{\mu(V)}-\mu(U)\right|\rightarrow0
$
as $\Delta\rightarrow\infty$.
Clearly $\psi$-mixing implies all the other kinds of mixing. The next strongest
mixing property is $\phi$-mixing, then comes strong mixing and uniform 
mixing is the weakest. The uniform strong mixing property is stronger that the 
strong mixing property but is implies by the dynamical $\phi$-mixing property
as we will see in Lemma~\ref{phi-mixing}. In fact if $\mu$ is strong mixing then 
the sets $S$ in Definition~\ref{definition} have to be of product form.

\vspace{3mm}

\noindent For a partition $\cal A$ we have the ($n$-th)
{\em information function} $I_n(x)=-\log\mu(A_n(x))$, where $A_n(x)$ denotes the
unique $n$-cylinder that contains the point $x\in \Omega$, whose moments are
$$
K_w({\cal A}) = \sum_{A\in{\cal A}} \mu(A)|\log\mu(A)|^w=\mathbb{E}(I_n^w),
$$
$w\ge0$ not necessarily integer. (For $w=1$ one traditionally writes 
$H({\cal A})=K_1({\cal A})=\sum_{A\in{\cal A}}-\mu(A)\log\mu(A)$.) If $\cal A$ is
finite then $K_w({\cal A})<\infty$ for all $w$. For infinite partitions
the theorem of Shannon-McMillan-Breiman requires that $H({\cal A})$ be finite~\cite{Car,Chu}.
In order to prove that the information function is lognormally distributed we will
require a larger than fourth moment  $K_w({\cal A})$ for some $w>4$ (not necessarily
integer) be finite.

\subsection{Results}

For $x\in\Omega$ we denote $A_n(x)$ the $n$-cylinder in ${\cal A}^n$
which contains the point $x$.
We are interested in the limiting behaviour of the distribution function
$$
\Xi_n(t)=\mu\left(\left\{x\in\Omega:
\frac{-\log\mu(A_n(x))-nh}{\sigma\sqrt{n}}\le t\right\}\right)
$$
for real valued $t$ and a suitable positive $\sigma$, where $h$ is the
metric entropy of $\mu$.
The Central Limit Theorem states that this quantity converges to
the normal distribution
$N(t)=\frac1{\sqrt{2\pi}}\int_{-\infty}^t e^{-s^2/2}\,ds$ as $n$ goes to
infinity if there exists a suitable $\sigma$ which is positive. Our main result is
the following theorem:

\begin{theorem}\label{CLT}
Let  $\mu$ be a uniformly strong mixing 
probability measure on $\Omega$ with respect to a countably finite, measurable and 
generating partition $\cal A$ which satisfies $K_w({\cal A})<\infty$ for some $w>4$. 
Assume that $\psi$  decays at
 least polynomially with power $>8+\frac{24}{w-4}$.

Then \\
{\bf (I)} The limit 
$$
 \sigma^2=\lim_{n\rightarrow\infty}\frac{K_2({\cal A}^n)-H^2({\cal A}^n)}n
 $$
exists and defines the variance of $\mu$. Moreover if the partition is infinite then
$\sigma$ is strictly positive.\\
{\bf (II)} If $\sigma>0$:
$$
\left|\Xi_n(t)-N(t)\right|\le C_0\frac1{n^\kappa}
$$
for all $t$ and all \\
(i) $\kappa<\frac1{10}-\frac35\frac{w}{(p+2)(w-2)+6}$ if $\psi$
 decays polynomially with power $p$,\\
(ii) $\kappa<\frac1{10}$ if $\psi$ decays faster than any power.
\end{theorem} 

\noindent The variance $\sigma^2$ is determined in Proposition~\ref{variance.convergence} 
and essentially only requires finiteness of the second moment $K_2({\cal A})$. 
In order to obtain the rate of convergence one usually needs a higher than second
moment of $I_n$. Since we use Stein's method we require the fourth moment be finite
(unlike in~\cite{HV2} where for finite partitions and $(\psi,f)$-mixing measures we only
needed bounds on the third moment).

Throughout the paper we shall assume that $K_w({\cal A})<\infty$ for some finite $w>4$. 
The case in which $w$ can be arbitrarily large (e.g.\ for finite partitions) is done
with minor modifications and yields the obvious result for the rate of convergence.
For simplicity's sake we assume in the proofs that the decay rate of $\psi$ is polynomial at
some finite power $p$. The case of hyper polynomial decay can be traced out with minor
modifications and yields the stated result. 

If the partition $\cal A$ is finite then $K_w({\cal A})<\infty$ 
for all $w$ and we obtain the following corollary:

\begin{corollary}\label{CLT.finite}
Let  $\mu$ be a uniformly strong mixing 
probability measure on $\Omega$ with respect to a finite, measurable and 
generating partition $\cal A$ and $\psi$  decays at
 least polynomially with power $>8+\frac{24}{w-4}$.

Then \\
{\bf (I)} The limit $\sigma^2=\lim_{n\rightarrow\infty}\frac1n(K_2({\cal A}^n)-H^2({\cal A}^n))$
exists (variance of $\mu$).\\
{\bf (II)} If $\sigma>0$: $\Xi_n(t)=N(t)+{\cal O}(n^{-\kappa})$
for all $t$ and
$\left\{\begin{array}{l}
\mbox{\it $\kappa<\frac1{10}-\frac3{5(p+2)}$
 if $\psi(\Delta)={\cal O}(\Delta^{-p})$, $\Delta\in\mathbb{N}$,}\\
\mbox{\it  $\kappa<\frac1{10}$ if $\psi$ decays hyper polynomially.}
\end{array}\right.$
\end{corollary}

\noindent By a result of Petrov~\cite{Pet} we now obtain the Law of the Iterated Logarithm from 
Theorem~\ref{CLT} by virtue of the error bound (better than  $\frac1{(\log n)^{1-\varepsilon}}$
(some $\varepsilon>0$) which are the ones required in \cite{Pet}).

\begin{corollary}
Under the assumptions of Theorem~\ref{CLT}:
$$
\limsup_{n\rightarrow\infty}\frac{I_n(x)-nh}{\sigma\sqrt{2n\log\log n}}=1
$$
almost everywhere.
\end{corollary}

\noindent A similar statement is true for the $\liminf$ where the limit is then equal to $-1$ 
almost everywhere.

\vspace{3mm}

\noindent Based on the Central Limit Theorem we also  get the weak invariance principle 
WIP (see section~4). Recently there has been a great 
interest in the WIP in relation to mixing properties of dynamical systems. 
For instance it has been obtained for a large class of observables and for a large
 class of dynamical systems by Chernov in \cite{ch}.  Other recent results are
 \cite{FMT,FHV,PS}. Those results however are typically for sums of sufficiently regular observables.
Here we prove the WIP for $I_n(x)$.

\begin{theorem}\label{wip}
Under the assumption of Theorem~\ref{CLT} 
 the information function $I_n$ satisfies the Weak
Invariance Principle (provided the variance $\sigma^2$ is positive).
\end{theorem}

\subsection{Examples}\label{examples}
{\bf (I) Bernoulli shift:} Let $\Sigma$ be the full shift space over the infinite alphabet $\mathbb{N}$ and
let $\mu$ be the Bernoulli measure generated by the positive weights $p_1,p_2,\dots$
($\sum_jp_j=1$). The entropy is then $h(\mu)=\sum_jp_j|\log p_j|$ and since
$K_2({\cal A})=\sum_ip_i\log^2p_i=\frac12\sum_{i,j}p_ip_j\left(\log^2 p_i +\log^2 p_j\right)$
we obtain that the variance is given by the following expression which is 
familiar from finite alphabet Bernoulli shifts:
$$
\sigma^2=K_2({\cal A})-h(\mu)^2=\frac12\sum_{i,j}p_ip_j\log^2\frac{p_i}{p_j}.
$$
We have used that the partition $\cal A$ is given by the cylinder sets whose first symbols
are fixed. Here we naturally assume that $\sum_ip_i\log^2p_i<\infty$.
If moreover $\sum_ip_i\log^4p_i<\infty$ then 
$$
\mathbb{P}\left(\frac{-\log\mu(A_n(x))-nh}{\sigma\sqrt{n}}\le t\right)=N(t)+{\cal O}(r^{-1/4})
$$
with exponent $\frac14$ which is a well known result for unbounded iid random variables.
With other techniques one can however weaken the moment requirement in this case. 

\vspace{3mm}

\noindent
{\bf (II) Markov shift:} Again let $\Sigma$ be the shift space over the infinite alphabet $\mathbb{N}$
and $\mu$ the Markov measure generated by an infinite probability vector $\vec{p}=(p_1,p_2,\dots)$
($p_j>0$, $\sum_jp_j=1$) and an infinite stochastic matrix $P$ ($\vec{p}P=\vec{p}$, $P\mathbf{1}=\mathbf{1}$).
The partition $\cal A$ is again the partition of single element cylinder sets. If 
$\vec{x}=x_1x_2\dots x_n$ is a word of length $n$ (we write $\vec{x}\in{\cal A}^n$)
then the measure of its cylinder set is 
$\mu(\vec{x})=p_{x_1}P_{x_1x_2}P_{x_2x_3}\cdots P_{x_{n-1}x_n}$. 
The metric entropy is $h(\mu)=\sum_{i,j}-p_iP_{ij}\log P_{ij}$~\cite{Wal} and
the variance~\cite{Rom,Yus} (see also Appendix) is
$$
\sigma^2=\frac12\sum_{ijk\ell}p_iP_{ij}p_kP_{k\ell}\log^2\frac{P_{ij}}{P_{k\ell}}
+4\sum_{k=2}^\infty\sum_{\vec{x}\in{\cal A}^k}\mu(\vec{x})\left(\log P_{x_1x_2}\log P_{x_{k-1}x_k}-h^2\right).
$$

\noindent
{\bf (III) Gibbs states:} The measure $\mu$ is a Gibbs state for the potential function $f$
if there exists a constant $c>1$ so that for every point $x\in \Omega$ and $n$ one has
$\mu(A_n(x))\in\left[\frac1c,c\right]e^{f^n(x)-nP(f)}$ where $P(f)$ is the pressure of $f$
and $f^n=f+f\circ T+\cdots+f\circ T^{n-1}$ is the $n$th ergodic sum of $f$. If $f$ is H\"older 
continuous and $T$ is and Axiom~A map, then $\mu$ is the unique equilibrium state.
In this case the CLT has been studied a great deal in particular for finite partitions since
standard techniques for sums of random variables can be applied (see e.g.~\cite{Bro,Gor,Liv,Pen}).
Note that (I) and (II) are special cases of Gibbs states.

\section{Variance and higher moments}

\subsection{Some basic properties}

\noindent Let us begin by showing that the uniform strong mixing property is implied
by the $\phi$-mixing property.

\begin{lemma} \label{phi-mixing}
$\phi$-mixing implies uniformly strong mixing.
\end{lemma}

\noindent {\bf Proof.} Let  $\mu$ be a left $\phi$-mixing probability measure
(the right $\phi$-mixing case is done in the same way).
That means, there exists a decreasing $\phi(\Delta)\rightarrow0$ as 
$\Delta\rightarrow\infty$ so that 
$$
\left|\left(\mu(B\cap C)-\mu(B)\mu(C)\right)\right|\le\phi(\Delta)\mu(B)
$$
for every $C$ in the $\sigma$-algebra generated by ${\cal C}=T^{-n-\Delta}{\cal A}^m$ and every
cylinder $B\in{\cal B}={\cal A}^n$ for all $n$ and $\Delta$. Let $S\subset{\cal B}\times{\cal C}$
and put $S_B$ for the interection of  $\{B\}\times{\cal C}$ with $S$. Then
$\left|\left(\mu(B\cap S_B)-\mu(B)\mu(S_B)\right)\right|\le\phi(\Delta)\mu(B)$
and
$$
\left|\sum_{(B,C)\in S}\left(\mu(B\cap C)-\mu(B)\mu(C)\right)\right|
\le\sum_{B\in{\cal B}}\left|\left(\mu(B\cap S_B)-\mu(B)\mu(S_B)\right)\right|
\le\sum_{B\in{\cal B}}\phi(\Delta)\mu(B)\le \phi(\Delta)
$$
implies that $\mu$ is uniformly strong mixing with $\psi=\phi$.\qed

\vspace{3mm}

\noindent The following estimate has previously been
shown for $\psi$-mixing measures (in which case they are exponential) in~\cite{GS} 
and for $\phi$-mixing measures in~\cite{Aba}.   
Denote by $A_n(x)$ the atom in ${\cal A}^n$ ($n=1,2,\dots$) which contains the point $x\in\Omega$.
(Abadi~\cite{Aba} also showed that in case~(II) the decay cannot in general be exponential.)

\begin{lemma}\label{cylinderestimate} Let $\mu$ be strong mixing.
Then there exists a constant $C_1$ so that for all $A\in{\cal A}_n$, $n=1,2,\dots$: \\
(I) $\mu(A)\le C_1n^{-p}$
if $\psi$ is polynomially decreasing with exponent $p>0$;\\
(II) $\mu(A)\le C_1\theta^{\sqrt{n}}$ for some $\theta\in(0,1)$
if $\psi$ is exponentially decreasing.
\end{lemma}

\noindent {\bf Proof.} Fix $m\ge1$ so that $a=\max_{A\in{\cal A}^m}\mu(A)$ is less than $\frac14$
and let $\Delta_1,\Delta_2,\dots$ be 
integers which will be determined below. We put $n_j=jm+\sum_{i=1}^{j-1}\Delta_i$
(put $\Delta_0=0$) and for $x\in \Omega$ let $B_j=A_m(T^{n_{j-1}+m}x)$ and 
put $C_k=\bigcap_{j=1}^kB_j$.
Then $A_{n_k}(x)\subset C_k$ and 
$$
\mu(C_{k+1})=\mu(C_k\cap B_{k+1})=\mu(C_k)\mu(B_{k+1})+\rho(C_k,B_{k+1})
$$
where the remainder term $\rho(C_k,B_{k+1})$ is by the mixing 
property in absolute value bounded by $\psi(\Delta_k)$.
Now we choose $\Delta_j$ so that $\psi(\Delta_j)\le a^{\frac{j}2+1}$.
Then $\mu(C_{k+1})\le\mu(C_k)a+a^{\frac{k}2+1}$ implies that
$\mu(C_k)\le c_0a^{\frac{k}2}$ (as $\sqrt{a}\le\frac12$) for some $c_0>0$.

\vspace{2mm}

\noindent  {\bf (I)} If $\psi$ decays polynomially with 
power $p$, i.e. $\psi(t)\le c_1t^{-p}$, then the condition $\psi(\Delta_j)\le a^{\frac{j}2+1}$
is satisfied if we put $\Delta_j=\left[c_2a^{-\frac{j}{2p}}\right]$ for a suitable constant $c_2>0$.
Consequently $n_k\le c_3a^{-\frac{k}{2p}}$ ($c_3\ge1$) and therefore
 $k\ge2p\frac{\log n_k}{|\log a|}$. Hence
 $$
 \mu(A_{n_k}(x))\le c_0a^{\frac{k}2}\le c_0a^{p\frac{\log n_k}{|\log a|}}
 \le c_4n_k^{-p}
 $$
 and from this one obtains $ \mu(A_n(x))\le c_5n^{-p}$ for all integers 
 $n$ (and some larger constant $c_5$).
 
\vspace{2mm}

\noindent  {\bf (II)}  If $\psi$ decays exponentially, i.e. $\psi(t)\le c_6\vartheta^t$
for some $\vartheta\in(0,1)$, then we choose 
$\Delta_j=\left[\frac{j}2\frac{\log a}{\log\vartheta}\right]$ and obtain 
$n_k\le mk+c_7k^2$, which gives us $k\ge c_8\sqrt{n_k}$ ($c_8>0$) and the
stretched exponential decay of the measure of cylindersets:
$$
\mu(A_n(x))\le c_9a^{c_8\sqrt{n}}.
$$
Now put $\theta=a^{c_8}$. \qed

\subsection{The information function and mixing properties}

\noindent The metric entropy $h$ for the invariant measure $\mu$ is
 $h=\lim_{n\rightarrow\infty}\frac1nH({\cal A}^n)$, where
  $\cal A$ is a generating partition of $\Omega$ (cf.\ \cite{Man}), provided $H({\cal A})<\infty$.
 For $w\ge1$ put $\eta_w(t)=t\log^w\frac1t$ ($\eta_w(0)=0$).
Then 
$$
K_w({\cal B}) = \sum_{B\in{\cal B}} \mu(B)|\log\mu(B)|^w
 = \sum_{B\in{\cal B}} \eta_w(\mu(B))
$$
for partitions $\cal B$.
Similarly one has the conditional quantity ($\cal C$  is a partition):
$$
K_w({\cal C}|{\cal B})
=\sum_{B\in{\cal B},C\in{\cal C}}\mu(B)\eta_w\left(\frac{\mu(B\cap C)}{\mu(B)}\right)
=\sum_{B,C}\mu(B\cap C)\left|\log\frac{\mu(B\cap C)}{\mu(B)}\right|^w.
$$

\begin{lemma}\label{general.subadditivity} \cite{HV2}
For any two partitions ${\cal B}, {\cal C}$ for which $K_w({\cal B}),K_w({\cal C})<\infty$
and $\mu(C)\le e^{-w}\;\;\forall\;\;C\in{\cal C}$:

\noindent (i) $K_w({\cal C}|{\cal B}) \leq K_w({\cal C})$,

\noindent (ii) $K_w({\cal B}\vee{\cal C})^{1/w}
\le K_w({\cal C}|{\cal B})^{1/w}+K_w({\cal B})^{1/w}$,

\noindent (iii) $K_w({\cal B}\vee{\cal C})^{1/w}
\le K_w({\cal C})^{1/w}+K_w({\cal B})^{1/w}.$
\end{lemma}

\noindent {\bf Proof.}
(i) Since $\eta_w(t)$ is convex and increasing on $[0,e^{-w}]$ and decreasing to zero
on $(e^{-w},1]$ we have $\sum_ix_i\eta_w(\alpha_i)\le\eta_w\left(\sum_ix_i\alpha_i\right)$
for weights $x_i\ge0$ ($\sum_ix_i=1$) and numbers $\alpha_i\in[0,1]$  which satisfy
$\sum_ix_i\alpha_i\le e^{-w}$. Hence 
$$
K_w({\cal C}|{\cal B})
=\sum_{B\in{\cal B},C\in{\cal C}}\mu(B)\eta_w\left(\frac{\mu(B\cap C)}{\mu(B)}\right)
\le\sum_C\eta_w\left(\sum_B\mu(B)\frac{\mu(B\cap C)}{\mu(B)}\right)
=\sum_C\eta_w\left(\mu(C)\right)
=K_w({\cal C}).
$$
\noindent (ii) The second statement follows from Minkowski's inequality on $L^w$-spaces:
\begin{eqnarray*}
K_w({\cal B}\vee{\cal C})^\frac1w&=&
\left(\sum_{B\in{\cal B},C\in{\cal C}}\mu(B\cap C)
\left|\log\mu(B\cap C)\right|^w\right)^\frac1w\\
&\leq&\left(\sum_{B,C}\mu(B\cap C)
\left|\log\frac{\mu(B\cap C)}{\mu(B)}\right|^w\right)^\frac1w
+\left(\sum_{B,C}\mu(B\cap C)\left|\log\mu(B)\right|^w\right)^\frac1w\\
&=&K_w({\cal C}|{\cal B})^\frac1w+K_w({\cal B})^\frac1w.
\end{eqnarray*}
(iii) This follows from (ii) and (i).
\hfill $\qed$

\begin{corollary}\label{K.estimate} Let $w\ge1$ and ${\cal A}$ so that $K_w({\cal A})<\infty$ and
 $\mu(A)\le e^{-w}\;\forall\;A\in{\cal A}$.
Then there exists a constant $C_2$ (depending on $w$) so that for all $n$
$$
K_w({\cal A}^n)\le C_2 n^w.
$$ 
\end{corollary}

\noindent {\bf Proof.} We want to use Lemma~\ref{general.subadditivity}(iii)
to show that the sequence $a_n=K_w({\cal A}^n)^{1/w}$, $n=1,2,\dots$,
 is subadditive.
The hypothesis of Lemma~\ref{general.subadditivity} is satisfied since
$\mu(A)\le e^{-w}$ for all $A\in{\cal A}$. 
 We thus obtain $K_w({\cal A}^{n+m})^\frac1w\le K_w({\cal A}^n)^\frac1w+K_w({\cal A}^m)^\frac1w$
 for all $n,m\ge1$ and therefore subadditivity of the sequence $a_n$.
 Since by assumption $K_w({\cal A})<\infty$ we get that
 the limit $\lim_{n\rightarrow\infty}\frac1nK_w({\cal A}^n)^{1/w}$ exists, is finite
and equals the $\inf$ (see e.g.~\cite{Wal}). 
\qed

\vspace{3mm}

\noindent 
The function $I_n$ has expected value $\mathbb{E}(I_n)=H({\cal A}^n)$, for which we also write
$H_n$, and variance $\sigma_n^2=\sigma^2(I_n)= K_2({\cal A}^n)-H_n^2$.
In general, if ${\cal B}$ is a partition then we write
$\sigma^2({\cal B})=K_2({\cal B})-H^2({\cal B})$ and similarly
for the conditional variance $\sigma^2({\cal C}|{\cal B})$. Let us define the 
function $J_{\cal B}$ by
  $J_{\cal B}(B)=-\log\mu(B)-H({\cal B})$ ($B\in{\cal B}$) then
  $\sigma^2({\cal B})=\sum_{B\in{\cal B}}\mu(B)J_{\cal B}(B)^2$
and $\int J_{\cal B}\,d\mu=0$. For two partitions $\cal B$ and 
$\cal C$ we put
 $$
 J_{{\cal C}|{\cal B}}(B\cap C)=\log\frac{\mu(B)}{\mu(B\cap C)}-H({\cal C}|{\cal B})
 $$
 for $(B,C)\in{\cal B}\times{\cal C}$. (This means 
 $J_{{\cal C}|{\cal B}}= J_{{\cal B}\vee{\cal C}}-J_{\cal B}$
 and $\sigma({\cal C}|{\cal B})=\sigma(J_{{\cal C}|{\cal B}})$.)
\begin{lemma}\label{variance.join}
Let $\cal B$ and $\cal C$ be two partitions.
Then
$$
\sigma({\cal B}\vee{\cal C})
\leq\sigma({\cal C}|{\cal B})+\sigma({\cal B}).
$$
\end{lemma}

\noindent {\bf Proof.} 
This follows from  Minkowski's inequality
$$
\sigma({\cal B}\vee{\cal C})
=\sqrt{\mu\left(J_{{\cal C}|{\cal B}}+J_{\cal B}\right)^2}
\le\sqrt{\mu(J_{{\cal C}|{\cal B}}^2)}+\sqrt{\mu(J_{\cal B}^2)}
=\sigma({\cal C}|{\cal B})+\sigma({\cal B}).
$$
\qed

\noindent As a consequence of Lemma~\ref{general.subadditivity}(i) one also has
$K_w({\cal B}\vee{\cal C}|{\cal B})=K_w({\cal C}|{\cal B})\le K_w({\cal C})$
which in particular implies 
$\sigma({\cal B}\vee{\cal C}|{\cal B})=\sigma({\cal C}|{\cal B})\le
\sqrt{K_2({\cal C})}$.
As before we put $\rho(B,C)=\mu(B\cap C)-\mu(B)\mu(C)$ (and in the following we often write
${\cal B}={\cal A}^n$ and ${\cal C}=T^{-\Delta-n}{\cal A}^n$ for integers $n, \Delta$).

The following technical lemma is central to get the variance of $\mu$ and bounds
on the higher moments of $J_n=I_n-H_n$.

\begin{lemma}\label{rho.remainder.estimate}
Let $\mu$ be uniformly strong mixing and assume that $K_w({\cal A})<\infty$
and $\mu(A)\le e^{-w}\;\forall\;A\in{\cal A}$ for some $w\ge1$. 
Then for every
$\beta>1$ and $a\in[0,w)$ there exists a constant $C_4$ so that 
$$
\sum_{B\in{\cal B}, C\in {\cal C}}\mu(B\cap C)\left|\log\left(1+\frac{\rho(B,C)}{\mu(B)\mu(C)}\right)\right|^a
\le C_4\left(\psi(\Delta)(m+n)^{(1+a)\beta}+(m+n)^{a\beta-w(\beta-1)}\right)
$$
for $\Delta<\min(n,m)$ and for all $n=1,2,\dots$.
(As before ${\cal B}={\cal A}^m$, ${\cal C}=T^{-\Delta-m}{\cal A}^n$.)
\end{lemma}

\noindent {\bf Proof.}  
Let $m, n$ and $\Delta$ be as in the statement and put 
$$
{\cal L}_\ell=\left\{(B,C)\in{\cal B}\times{\cal C}:
\;\; 2^{\ell-1}<1+\frac{\rho(B,C)}{\mu(B)\mu(C)}\le2^\ell\right\}
$$
$\ell\in\mbox{\bf Z}$.
Using the  strong mixing property we obtain
$$
\sum_{B\in{\cal B}, C\in {\cal C}}\mu(B\cap C)\left|\log\left(1+\frac{\rho(B,C)}{\mu(B)\mu(C)}\right)\right|^a
=\sum_{\ell=-\infty}^\infty L_\ell(|\ell|+{\cal O}(1))^a
$$
where $L_\ell=\sum_{(B,C)\in{\cal L}_\ell}\mu(B\cap C)$.
Since $\rho(B,C)={\cal O}(1)(2^\ell-1)\mu(B)\mu(C)$ we get
${\cal O}(\psi(\Delta))=\sum_{(B,C)\in{\cal L}_\ell}\rho(B,C)={\cal O}(1)(2^\ell-1)L_\ell^\times$
where $L_\ell^\times=\sum_{(B,C)\in{\cal L}_\ell}\mu(B)\mu(C)$. Hence
for $\ell>0$ one obtains $L_\ell^\times={\cal O}(\psi(\Delta))2^{-\ell}$ and
if $\ell<0$ then $L_\ell^\times={\cal O}(\psi(\Delta))$.
Also note that if $\ell=0$ then 
$\left|\log\left(1+\frac{\rho(B,C)}{\mu(B)\mu(C)}\right)\right|
={\cal O}\left(\frac{\rho(B,C)}{\mu(B)\mu(C)}\right)$ and
$$
\sum_{(B,C)\in{\cal L}_0}\mu(B\cap C)\left|\log\left(1+\frac{\rho(B,C)}{\mu(B)\mu(C)}\right)\right|^a
={\cal O}(1)\sum_{(B,C)\in{\cal L}_0}\rho(B,C)={\cal O}(\psi(\Delta)).
$$
We separately estimate (i) for $\ell\ge1$ and (ii) for $\ell\le-1$:\\
{\bf (i)} Since $\mu(B\cap C)=\left(1+\frac{\rho(B,C)}{\mu(B)\mu(C)}\right)\mu(B)\mu(C)$ we get
for $\ell \ge1$:
$$
2^{\ell-1}L_\ell^\times=\sum_{(B,C)\in{\cal L}_\ell}\mu(B)\mu(C)2^{\ell-1}
\le L_\ell\le\sum_{(B,C)\in{\cal L}_\ell}\mu(B)\mu(C)2^\ell=2^\ell L_\ell^\times
$$
Thus
$$
\sum_{\ell=1}^{(m+n)^\beta}\ell^a L_\ell
\le\sum_{\ell=1}^{(m+n)^\beta}\ell^a2^\ell L_\ell^\times
\le\sum_{\ell=1}^{(m+n)^\beta}\ell^a\frac{2^\ell}{2^\ell-1}\psi(\Delta)
\le c_1\psi(\Delta)(m+n)^{(1+a)\beta}.
$$
For $\ell>(m+n)^\beta$ we use that $\mu(B\cap C)\ge2^{\ell-1}\mu(B)\mu(C)$ on ${\cal L}_\ell$
which implies $\mu(B)\mu(C)\le2^{1-\ell}$ and 
$\mu(B\cap C)\le\min(\mu(B),\mu(C))\le2^{-\frac{\ell-1}2}$. Hence, on  ${\cal L}_\ell$
one has $|\log\mu(B\cap C)|\ge(\ell-1) \log\sqrt{2}$. Similarly to the previous lemma put
$$
D_k=\bigcup_{(B,C)\in{\cal B}\times{\cal C},\;k-1<|\log\mu(B\cap C)|\le k}(B\cap C)
$$
and use Corollary~\ref{K.estimate} to get
(as $K_w({\cal A}^{n+m+\Delta})\ge K_w({\cal B}\vee{\cal C})$)
$$
C_2(w)(n+m+\Delta)^w\ge K_w({\cal B}\vee{\cal C})
\ge\sum_{k=1}^\infty\mu(D_k)(k-1)^w
\ge c_2(n+m)^{\beta (w-a)}\sum_{k=[(n+m)^\beta]+1}^\infty\mu(D_k)k^{a}.
$$
We thus obtain (using that $\Delta<\min(n,m)$)
\begin{eqnarray*}
\sum_{\ell=(n+m)^\beta}^\infty \ell^a L_\ell
&\le&\frac1{\log\sqrt{2}}\sum_{(B,C)\in{\cal B}\times{\cal C},\;|\log\mu(B\cap C)|\ge(n+m)^\beta}
 \left|\log\mu(B\cap C)\right|^a\mu(B\cap C)\\
&\le&\frac1{\log\sqrt{2}}\sum_{k=[(n+m)^\beta]+1}^\infty k^a\mu(D_k)\\
&\le&c_3\frac{(n+m)^{a\beta}}{(n+m)^{(\beta-1) w}}
\end{eqnarray*}
for some $c_3$ (which depends on $w$).\\
{\bf (ii)} For negative values of $\ell$ we use 
$L_\ell\le2^\ell L_\ell^\times \le c_42^\ell\psi(\Delta)$ which gives
$$
\sum_{\ell=-\infty}^0 |\ell|^a L_\ell
\le c_4\sum_{\ell=0}^\infty \ell^a2^{-\ell}\psi(\Delta)
\le c_5\psi(\Delta).
$$
Combining (i) and (ii) yields 
$$
\sum_{\ell=-\infty}^\infty L_\ell(|\ell|+{\cal O}(1))^a
\le (c_1+c_4)\psi(\Delta)(m+n)^{(1+a)\beta}+c_3(m+n)^{a\beta-w(\beta-1)}
$$
which concludes the proof.\qed

\subsection{Entropy}
The main purpose of this section is to obtain rates of convergence for the 
entropy  (Lemma~\ref{entropy.approximation}).

\begin{lemma}\label{entropy.additivity}
Under the assumptions of Lemma~\ref{rho.remainder.estimate}
for every $\beta>1$ there exists a constant $C_5$ so that for all $n$:
$$
\left|H({\cal B}\vee{\cal C})-(H({\cal B})+H({\cal C}))\right|
\le C_5\left(\psi(\Delta)n^{2\beta}+n^{\beta-(\beta-1)w}\right),
$$
where ${\cal B}={\cal A}^n$, ${\cal C}=T^{-\Delta-n}{\cal A}^n$.
\end{lemma}

\noindent {\bf Proof.} 
Using the uniform strong mixing property $\mu(B\cap C)=\mu(B)\mu(C)+\rho(B,C)$ we obtain
\begin{eqnarray*}
H({\cal B}\vee{\cal C})&=&\sum_{B\in{\cal B}, C\in {\cal C}} \mu(B\cap
C) \log \frac{1}{\mu (B\cap C)}\\
&=&\sum_{B,C}\mu(B\cap C)
\left(\log \frac1{\mu(B)} + \log \frac1{\mu(C)} -\log\left(1+\frac{\rho(B,C)}{\mu(B)\mu(C)}\right)\right)\\
&=& H({\cal B}) + H({\cal C})  + E,
\end{eqnarray*}
where by Lemma~\ref{rho.remainder.estimate} (with $a=1$) 
$$
E=-\sum_{B\in{\cal B}, C\in {\cal C}}\mu(B\cap C)\log\left(1+\frac{\rho(B,C)}{\mu(B)\mu(C)}\right)
={\cal O}\left(\psi(\Delta)n^{2\beta}+n^{\beta-(\beta-1)w}\right).
$$
This proves the lemma.\qed

\begin{lemma} \label{entropy.approximation}
Under the assumptions of Lemma~\ref{rho.remainder.estimate}
 there exists a constant $C_6$ so that ($H_m=H({\cal A}^m)$)
$$
\left|\frac{H_m}m-h\right|\le C_6\frac1{m^\gamma}
$$
for all $m$,
where $\gamma\in(0,1-\frac{2w}{p(w-1)})$ if $\psi$ decays polynomially with power 
$p>\frac{2w}{w-1}$ and $\gamma\in(0,1)$ if $\psi$ decays faster than polynomially. 
\end{lemma}

\noindent {\bf Proof.}  Let $m$ be an integer. Let  ${\cal B}={\cal A}^{u-\Delta}$,
${\cal C}=T^{-u}{\cal A}^{u-\Delta}$ and ${\cal D}=T^{-u}{\cal A}^{2\Delta}$, then 
by Lemma~\ref{entropy.additivity}:
$$
H_{2u}
=2H_{u-\Delta}+{\cal O}(H_{2\Delta})+{\cal O}\left(\psi(\Delta)u^{2\beta}+n^{\beta-(\beta-1)w}\right).
$$
If we choose $\delta\in(\frac{2w}{p(w-1)},1)$  and put $\beta=\frac{w}{w-1}$ then 
$\Delta={\cal O}(u^\delta)$ implies that $\psi(\Delta)u^{2\beta}+u^{\beta-(\beta-1)w}={\cal O}(1)$.
With $\Delta=[u^\delta]$ we thus obtain $H_{2u}=2H_u+{\cal O}(\Delta)=2H_u+{\cal O}(u^\delta)$
as  $H_{2\Delta}={\cal O}(\Delta)$ and $H_{u-\Delta}=H_u+{\cal O}(\Delta)$.
Iterating this estimate yields the following bound along exponential progression:
$$
H_{2^im}
=2^iH_m+\sum_{j=0}^{i-1}2^{i-1-j}{\cal O}\left((2^jm)^\delta\right)
=2^iH_m+{\cal O}\left(m^\delta2^i\right).
$$
To get bounds for arbitrary (large) integers $n$ we do the following dyadic argument: 
Let $n=km+r$ where 
$0\le r<m$ and consider the binary expansion of: 
$k=\sum_{i=0}^\ell\epsilon_i2^i$, where
$\epsilon_i=0,1$ ($\epsilon_\ell=1$, $\ell=[\log_2k]$). We also put
$k_j=\sum_{i=0}^{j}\epsilon_i2^i$ ($k_\ell=k$). Obviously
$k_j=k_{j-1}+\epsilon_j2^j\le2^{j+1}$.
If $\epsilon_j=1$ then we separate the `first' block of length $k_{j-1}m$ from the
`second' block of length $2^jm$ by a gap of
length $2[(k_{j-1}m)^\delta]$ which we cut away in equal parts from 
the two adjacent blocks). We thus obtain ($H_0=0$)
$$
H_{mk_j}=H_{\epsilon_j2^jm+k_{j-1}m}
=H_{\epsilon_j2^jm}+H_{k_{j-1}m}+{\cal O}(\epsilon_j(k_{j-1}m)^\delta)
=H_{\epsilon_j2^jm}+H_{k_{j-1}m}+{\cal O}(\epsilon_j(2^jm)^\delta)
$$
for $j=0,1,\dots,\ell-1$.
Iterating this formula and summing over $j$ yields
$$
H_{km}=\sum_{j=0}^\ell\epsilon_j\left(2^jH_m
+{\cal O}\left(m^\delta2^j\right)\right)
=kH_m+{\cal O}\left(m^\delta 2^\ell\right).
$$
The contribution made by the remainder of length $r$ is easily bounded by
$$
\left|H_n-H_{km}\right|
\le\sigma\left({\cal A}^n|{\cal A}^{km}\right)
\le c_1r\le c_1m.
$$
Consequently
$$
H_n=kH_m+{\cal O}\left(m^\delta 2^\ell\right)+{\cal O}(m)=kH_m+{\cal O}\left(m^\delta k\right)
$$
as $2^\ell\le k\le2^{\ell+1}$. 
Dividing by $n$ and letting $n$ go to infinity ($k\rightarrow\infty$) yields 
$$
h=\liminf_{n\rightarrow\infty}\frac{H_n}n
=\frac{H_m}m+{\cal O}\left(m^{\delta-1}\right)
$$
for all $m$ large enough.\qed

\subsection{The variance}
In this  section we prove part (I) of Theorem~\ref{CLT} and moreover obtain
convergence rates which will be needed to prove part (II) in section~4.

\begin{proposition}\label{variance.convergence}
Let $\mu$ be uniformly strong mixing and assume that $K_w({\cal A})<\infty$
and $\mu(A)\le e^{-w}\;\forall\;A\in{\cal A}$ for some $w>2$. 
Assume that $\psi$ is at least polynomially decaying with power $p>6+\frac{8}{w-2}$.
Then the limit
$$
\sigma^2=\lim_{n\rightarrow\infty}\frac1n\sigma^2({\cal A}^n)
$$
exists and is finite. Moreover for every $\eta<\eta_0=2\frac{(p-2)(w-2)}{(w-2)(p+2)+8}$
there exists a constant $C_7$ so that for all $n\in\mathbb{N}$:
$$
\left|\sigma^2-\frac{\sigma^2({\cal A}^n)}n\right|\le \frac{C_7}{n^\eta}.
$$
Moreover, if the partition $\cal A$ is infinite, then $\sigma$ is strictly positive.
\end{proposition}

\noindent {\bf Proof.}  With ${\cal B}={\cal A}^n, {\cal C}=T^{-n-\Delta}{\cal A}^n$ we have
 by Lemma~\ref{entropy.additivity}
$
H({\cal B}\vee{\cal C})
=H({\cal B})+H({\cal C})+{\cal O}\left(\psi(\Delta)n^{2\beta}+n^{\beta-(\beta-1)w}\right),
$
and get for the variance
\begin{eqnarray*}
\sigma^2({\cal B}\vee{\cal C})&=&\sum_{B\in{\cal B},C\in{\cal C}}
\mu(B\cap C)\left(\log\frac1{\mu(B\cap C)}
-H({\cal B}\vee{\cal C})\right)^2\\
&=&\sum_{B,C}
\mu(B\cap C)\left(J_{\cal B}(B)+J_{\cal C}(C)
+{\cal O}\left(\psi(\Delta)n^{2\beta}+n^{\beta-(\beta-1)w}\right)
-\log\left(1+\frac{\rho(B,C)}{\mu(B)\mu(C)}\right)\right)^2.
\end{eqnarray*}
By Minkowski's inequality:
$$
\left|\sigma({\cal B}\vee{\cal C})-\sqrt{ E({\cal B},{\cal C})} \right|
\le c_1\left(\psi(\Delta)n^{2\beta}+n^{\beta-(\beta-1)w}\right)+\sqrt{F({\cal B},{\cal C})}
$$
($c_1>0$) where (by Lemma~\ref{rho.remainder.estimate} with $a=2$)
$$
F({\cal B},{\cal C})= \sum_{B\in{\cal B},C\in{\cal C}}
\mu(B\cap C)\log^2\left(1+\frac{\rho(B,C)}{\mu(B)\mu(C)}\right)
\le c_2\left(\psi(\Delta)n^{3\beta}+n^{2\beta-(\beta-1)w}\right),
$$
and
\begin{eqnarray*}
E({\cal B},{\cal C})&=& \sum_{B\in{\cal B},C\in{\cal C}}
\mu(B\cap C)\left(J_{\cal B}(B)+J_{\cal C}(C)\right)^2\\
&=&\sum_{B,C}\mu(B\cap C)
\left(J_{\cal B}(B)^2+J_{\cal C}(C)^2\right)+2G({\cal B},{\cal C}) \\
&=& \sigma^2({\cal B})+ \sigma^2({\cal C})+ 2G({\cal B},{\cal C}).
\end{eqnarray*}
Since $J_{\cal B}$ and $J_{\cal C}$ have average zero the remainder term
\begin{eqnarray*}
G({\cal B},{\cal C})&=&\sum_{B\in{\cal B},C\in{\cal C}}
\mu(B\cap C)J_{\cal B}(B)J_{\cal C}(C)\\
&=&\sum_{B,C}(\mu(B)\mu(C)+\rho(B,C))J_{\cal B}(B)J_{\cal C}(C)\\
&=&\sum_{B,C}\rho(B,C)J_{\cal B}(B)J_{\cal C}(C)
\end{eqnarray*}
which is estimated using Schwarz' inequality as follows
$$
|G({\cal B},{\cal C})|\le\sum_{B,C}|\rho(B,C)|\cdot|J_{\cal B}(B)|\cdot|J_{\cal C}(C)|\\
\le\psi(\Delta)\sigma({\cal B})\sigma({\cal C}).
$$
Hence
\begin{equation}\label{sigma.addition}
\sigma ({\cal B}\vee{\cal C})
\leq\sqrt{\sigma^2({\cal C})+ \sigma^2({\cal B})
+\psi(\Delta)\sigma({\cal B})\sigma({\cal C})}
+ c_4\sqrt{\psi(\Delta)n^{3\beta}+n^{2\beta-(\beta-1)w}}.
\end{equation}
Next we fill the gap of length $\Delta$ for which we use 
Lemma~\ref{variance.join} and Corollary~ \ref{K.estimate}
$$
|\sigma({\A}^{2n+\Delta})- \sigma({\B}\vee{\C})|
\le\sigma(T^{-n}{\A}^\Delta|{\B}\vee{\C})
\le \sqrt{K_2(T^{-n}{\A}^{\Delta})}
=\sqrt{K_2({\A}^{\Delta})}
\le c_5 \Delta.
$$
Since by assumption $\psi(\Delta)\le c_6\Delta^{-p}$ for some $p>6+\frac{8}{w-2}$ we take can
$\delta=\frac{4w}{(p+2)(w-2)+8}$ and $\beta=\frac{2+p}4\delta$ (in particular
$\delta<\frac12$). Then, with $\Delta=[n^\delta]$ we get 
$\psi(\Delta)n^{4\beta}+n^{2\beta-(\beta-1)w}\le \Delta^2$. 
Therefore, as $\sigma({\cal B})= \sigma({\C})=\sigma_n$ (where $\sigma_n=\sigma({\cal A}^n)$),
one has
$$
\sigma_{2n+[n^\delta]}\le \sqrt{(2+\psi(\Delta))\sigma^2_n+c_7n^{2\delta}}
\le \sqrt{2\sigma^2_n+c_7n^{2\delta}},
$$
where in the last step we took advantage of the a priori estimates from Corollary~\ref{K.estimate}
$\sigma^2({\cal A}^n)\le K_2({\cal A}^n)\le C_2n^2$ and the choice of $\delta$ which implies
that $\psi(\Delta)n^2={\cal O}(1)$.
Since $2\delta<1$ one has
$\sigma^2_k\le c_8k$ for all $k$ and some constant $c_8$.
Given $n_0$ let us put recursively $n_{j+1}=2n_j+[n_j^\delta]$ ($j=0,1,2,\dots$).
Then $2^jn_0\le n_j\le2^jn_0\prod_{i=0}^{j-1}\left(1+\frac12n_i^{\delta-1}\right)$
where the product is bounded by
$$
\prod_{i=0}^{j-1}\left(1+\frac12n_i^{\delta-1}\right)
\le\prod_{i=0}^{j-1}\left(1+\frac{1}{n_0^{1-\delta}2^{(1-\delta)i+1}}\right)
\le\exp\frac{c_9}{n_0^{1-\delta}}.
$$
In the same fashion one shows that 
$\left|\sigma_{n_{j+1}}^2-2\sigma^2_{n_j}\right|\le c_7n_j^{2\delta}$
implies 
$$
2^j\sigma^2_{n_0}\exp-\frac{c_{10}}{n_0^{1-2\delta}}
\le\sigma_{n_j}^2
\le 2^j\sigma^2_{n_0}\exp\frac{c_{10}}{n_0^{1-2\delta}}.
$$
Hence
$$
\frac{2^j\sigma^2_{n_0}}{2^jn_0}\exp-\left(\frac{c_{10}}{n_0^{1-2\delta}}+\frac{c_9}{n_0^{1-\delta}}\right)
\le\frac{\sigma^2_{n_j}}{n_j}
\le\frac{2^j\sigma^2_{n_0}}{2^jn_0}\exp\frac{c_{10}}{n_0^{1-2\delta}},
$$
which simplifies to 
\begin{equation}\label{variance.linearity}
\frac{\sigma^2_{n_j}}{n_j}
=\frac{\sigma^2_{n_0}}{n_0}\left(1+{\cal O}\left(\frac1{n_0^{1-2\delta}}\right)\right)
=\frac{\sigma^2_{n_0}}{n_0}+{\cal O}\left(\frac1{n_0^{2-2\delta}}\right).
\end{equation}
As $w>2$ one has $\sigma_{n_0}<\infty$.
Taking $\limsup$ as $j\rightarrow\infty$ and $n_0\rightarrow\infty$ shows that the limit 
$\sigma^2=\lim_n\frac{\sigma^2_n}n$ exists and satisfies moreover 
$\left|\sigma^2-\frac{\sigma^2_n}n\right|\le C_7n^{-(2-2\delta)}$
for some $C_7$. Now we obtain the statement in the proposition for all
 $\eta<2-2\delta=2\frac{(p-2)(w-2)}{(p+2)(w-2)+8}$.
 
 In order to prove the last statement of the proposition let $\cal A$ be an infinite partition.
 If we choose $n_0$ large enough so that the error term ${\cal O}(n_0^{-(1-2\delta)})$ 
 in equation~(\ref{variance.linearity}) is $<\frac12$,
  then $\sigma^2_{n_j}>\frac12n_j\sigma^2_{n_0}$ for all $j$. 
 Since 
\begin{eqnarray*}
\sigma_{n_0}^2&=&\sum_{A\in{\cal A}^{n_0}} \mu(A)\log^2\mu(A)
-\sum_{A,B\in{\cal A}^{n_0}} \mu(A)\mu(B)\log\mu(A)\log\mu(B)\\
&=&\frac12\sum_{A,B} \mu(A)\mu(B)(\log^2\mu(A)+\log^2\mu(B))
-\sum_{A,B} \mu(A)\mu(B)\log\mu(A)\log\mu(B)\\
&=&\frac12\sum_{A,B\in{\cal A}^{n_0}} \mu(A)\mu(B)\log^2\frac{\mu(A)}{\mu(B)}
\end{eqnarray*}
we conclude that $\sigma_{n_0}^2>0$. Hence $\sigma^2=\lim_n\frac{\sigma^2_n}n$ 
is strictly positive.
\qed

\vspace{3mm}

\noindent {\bf Remarks:} (i) It is well known that for finite partitions the measure has variance zero
if it is a Gibbs state for a potential which is a coboundary.  \\
(ii) This proposition implies in particular that the limit
$\lim_{n\rightarrow\infty}\frac1{n^2}K_2({\cal A}^n)$
exists and is equal to $h^2$. \\
(iii) An application of Chebycheff's inequality gives
the large deviation type estimate ($\sigma_n=\sigma(J_n)$)
$$
\mathbb{P}\left(\frac1nJ_n(x)\ge t\right)
\le\frac{\sigma^2_n}{n^2t^2}= {\cal O}\left(\frac1{nt^2}\right).
$$

\subsection{Higher order moments}
\noindent In the proof of Theorem \ref{CLT} part~(II) we will need estimates
on the  third and fourth moments of $J_n$.  We first estimate the fourth moment 
and then use H\"older's inequality to bound the third moment. Denote by 
$$
M_w({\cal B})=\sum_{B\in{\cal B}}\mu(B)|J_{\cal B}(B)|^w.
$$
the $w$th (absolute) moment of the function $J_{\cal B}$.
By Minkowski's inequality
$$
M_4^\frac14({\cal B}\vee{\cal C})
=\sqrt[4]{\mu\left(J_{{\cal C}|{\cal B}}+J_{\cal B}\right)^4}
\le\sqrt[4]{\mu(J_{{\cal C}|{\cal B}}^4)}+\sqrt[4]{\mu(J_{\cal B}^4)}
=M_4^\frac14({\cal C}|{\cal B})+M_4^\frac14({\cal B}),
$$
where $M_w({\cal C}|{\cal B})=\sum_{B\in{\cal B},C\in{\cal C}}\mu(B\cap C)|J_{{\cal C}|{\cal B}}(B\cap C)|^w$
are the conditional moments.
It follows from Corollary~\ref{K.estimate} that the absolute moments for the
joins ${\cal A}^n$ can roughly be bounded by $M_w({\cal A}^n)\le K_w({\cal A}^n)\le C_2n^w$.
This estimate however is useless to prove Theorem~\ref{CLT} and the purpose of the 
next proposition is to reduce the exponent $w$ to $\frac12w$ 
in the cases $w=3,4$. One can of course get these improved estimates also for
$w$ larger than $4$ (as long as $K_w({\cal A})<\infty$) but we don't need those higher order moments here.

\begin{proposition}\label{fourth.moment}
Let $\mu$ be uniformly strong mixing and assume that $K_w({\cal A})<\infty$
and $\mu(A)\le e^{-w}\;\forall\;A\in{\cal A}$ for some $w>4$. 
Also assume that $\psi$ decays at least polynomially with power 
$>8+\frac{24}{w-4}$.
Then there exists a constant $C_8$ so that for all $n$
$$
M_4({\cal A}^n)\le C_8n^2$$
\end{proposition}

\noindent {\bf Proof.} 
With ${\cal B}={\cal A}^n$, ${\cal C}=T^{-\Delta-n}{\cal A}^n$ we get
(by Lemma~\ref{entropy.additivity})
 $H({\cal B}\vee{\cal C})=H({\cal B}) + H({\cal C})  + {\cal O}(\psi(\Delta)n^{2\beta}+n^{1-(\beta-1)w})$
and with Minkowsky's inequality (on $L^4$ spaces)
\begin{eqnarray*}
M_4^\frac14({\cal B}\vee{\cal C})&=&\left(\sum_{B\in{\cal B},C\in{\cal C}}
\mu(B\cap C)\left(\log\frac1{\mu(B\cap C)}
-H({\cal B}\vee{\cal C})\right)^4\right)^\frac14\\
&\leq&E_4^\frac14({\cal B},{\cal C})+{\cal O}\left(\psi(\Delta)n^{2\beta}+n^{\beta-(\beta-1)w}\right)
+F_4^{\frac14}({\cal B},{\cal C}),
\end{eqnarray*}
where by Lemma~\ref{rho.remainder.estimate} (with $a=4$)
$$
F_4({\cal B},{\cal C})= \sum_{B\in{\cal B},C\in{\cal C}}
\mu(B\cap C)\log^4\left(1+\frac{\rho(B,C)}{\mu(B)\mu(C)}\right)
= {\cal O}\left(\psi(\Delta)n^{5\beta}+n^{4\beta-(\beta-1)w}\right)
$$and
\begin{eqnarray*}
E_4({\cal B},{\cal C})&=& \sum_{B\in{\cal B},C\in{\cal C}}
\mu(B\cap C)\left(J_{\cal B}(B)+J_{\cal C}(C)\right)^4\\
&=&M_4({\cal B})+M_4({\cal C}) 
+\sum_{B,C}\mu(B\cap C)
\left(4J_{\cal B}(B)^3J_{\cal C}(C)+6J_{\cal B}(B)^2J_{\cal C}(C)^2+4J_{\cal B}(B)J_{\cal C}(C)^3\right).
\end{eqnarray*}
We look individually at the terms in the bracket:
\begin{eqnarray*}
\left|\sum_{B\in{\cal B},C\in{\cal C}}\mu(B\cap C)J_{\cal B}(B)^3J_{\cal C}(C)\right|
&=&\left|\sum_{B,C}(\mu(B)\mu(C)+\rho(B,C))J_{\cal B}(B)^3J_{\cal C}(C)\right|\\
&\le&\sum_{B,C}|\rho(B,C)|\cdot|J_{\cal B}(B)|^3|J_{\cal C}(C)|\\
&\le&\psi(\Delta)M_3({\cal B})\sigma({\cal C}).
\end{eqnarray*}
because $J_{\cal B}$ and $J_{\cal C}$ have zero average Schwarz inequality. In the same way we get
$$
\left|\sum_{B\in{\cal B},C\in{\cal C}}\mu(B\cap C)J_{\cal B}(B)J_{\cal C}(C)^3\right|
\le\psi(\Delta)\sigma({\cal B})M_3({\cal C}).
$$
Moreover
\begin{eqnarray*}
\sum_{B\in{\cal B},C\in{\cal C}}\mu(B\cap C)J_{\cal B}(B)^2J_{\cal C}(C)^2
&=&\sum_{B,C}(\mu(B)\mu(C)+\rho(B,C))J_{\cal B}(B)^2J_{\cal C}(C)^2\\
&=&\sigma^2({\cal B})\sigma^2({\cal C})+G({\cal B},{\cal C}),
\end{eqnarray*}
where
$$
|G({\cal B},{\cal C})|
=\left|\sum_{B\in{\cal B},C\in{\cal C}}\rho(B,C)J_{\cal B}(B)^2J_{\cal C}(C)^2\right|
\le\psi(\Delta)\sigma^2({\cal B})\sigma^2({\cal C})
$$
Thus 
$$
E_4({\cal B},{\cal C})=M_4({\cal B})+M_4({\cal C})+(6+\psi(\Delta))\sigma^2({\cal B})\sigma^2({\cal C})
+\psi(\Delta)\left(M_3({\cal B})\sigma({\cal C})+\sigma({\cal B})M_3({\cal C})\right).
$$

As  $\sigma^2({\cal B})=\sigma^2({\cal C})=\sigma^2_n\le c_1 n$
(Proposition~\ref{variance.convergence}) and since by assumption $\psi(\Delta)={\cal O}(\Delta^{-p})$
where $p>8+\frac{24}{w-4}$ we can choose $\beta=1+\frac2{w-4}$, $\delta=\frac1p(4+\frac{12}{w-4})$
and put $\Delta=[n^\delta]$. This implies $\Delta<\sqrt{n}$ (as $\delta<\frac12$) and
$\psi(\Delta)n^{6\beta}+n^{4\beta-(\beta-1)w}={\cal O}(n^2)$.
Using the a priori estimates $M_3({\cal A}^n)\le K_3({\cal A}^n)\le C_2n^3$ we obtain in particular that
$\psi(\Delta)\left(M_3({\cal B})\sigma({\cal C})+\sigma({\cal B})M_3({\cal C})\right)
={\cal O}(n^2)$ and therefore
$$
M_4^\frac14({\cal B}\vee{\cal C})
=\sqrt[4]{M_4({\cal C})+M_4({\cal B})+c_2n^2}
+ {\cal O}\left(\psi(\Delta))n^{2\beta}+n^{\beta-(\beta-1)w}\right),
$$
where the error term on the right hand side is ${\cal O}(n^{-3})$.
To fill in the gap of length $\Delta$ we use Lemma \ref{variance.join}
and the estimate on $K_4$ (Corollary~\ref{K.estimate}):
$$
\left|M_4^\frac14({\A}^{2n+\Delta})-M_4^\frac14({\B}\vee{\C})\right|\le
M_4^\frac14({\A}^{2n+\Delta} |{\B}\vee{\C})
\le K_4^\frac14({\A}^\Delta)
\le c_3\Delta.
$$
Hence
$$
M_4^\frac14({\A}^{n'})\le \sqrt[4]{2M_4({\cal A}^n)+c_2n^2}+c_3\Delta
\le\sqrt[4]{2M_4({\cal A}^n)+c_4n^2}
$$
(as $\Delta\le\sqrt{n}$), and by
induction  $M_4({\A}^k)\le C_8k^2$ (with $C_8\ge c_4/2$).
\qed

\vspace{3mm}

\noindent A H\"older estimate lets us now estimate the third 
absolute moments of $J_n$ as follows.

\begin{corollary}\label{third.moment}
Under the assumptions of Proposition~\ref{fourth.moment}
there exists a constant $C_9$ so that  for all $n$
$$
M_3({\cal A}^n)\le C_9n^{\frac32}.
$$
\end{corollary}

\section{Proof of Theorem~\ref{CLT} (CLT for Shannon-McMillan-Breiman)}

As before $N(t)$ denotes the normal distribution with zero mean and variance one.
We will first show the following result (in which $nh$ has been replaced by $H_n$
and $\sigma\sqrt{n}$ by $\sigma_n$).

\begin{theorem}\label{CLTn}
Under the assumptions of Theorem~\ref{CLT} one has:\\
{\bf (I)} The limit $\sigma^2=\lim_{n\rightarrow\infty}\frac1n(K_2({\cal A}^n)-H_n^2)$ exists
(and is positive if $|{\cal A}|=\infty$).\\
{\bf (II)} If $\sigma>0$ then
$$
\mathbb{P}\left(\frac{I_n-H_n}{\sigma_n}\le t\right)=N(t)+{\cal O}\left(\frac1{n^\kappa}\right)
$$
for all $t$ and all \\
(i) $\kappa<\frac1{10}-\frac35\frac{w}{(p+2)(w-2)+6}$ if $\psi$
 decays polynomially with power $p$,\\
(ii) $\kappa<\frac1{10}$ if $\psi$ decays hyper polynomially.
\end{theorem}
 
\noindent  {\bf Proof of Theorem~\ref{CLTn}.} 
It is enough to prove the theorem with the partition $\cal A$ replaced
by one of its joins ${\cal A}^k$ for some $k$. Since by Lemma~\ref{cylinderestimate}
$\mu(A)\le e^{-w}\;\forall \;A\in{\cal A}^k$ for some $k\ge1$ we therefore 
replace the original partition by ${\cal A}^k$ and will henceforth assume
that $\mu(A)\le e^{-w}$ for all $A\in{\cal A}$.

Theorem~\ref{CLTn} part~(I) follows from 
Proposition~\ref{variance.convergence}. For the proof of part~(II) let us 
assume that $\sigma$ is positive.
 We will use Stein's method to 
prove the CLT in the form of the following proposition which is modelled after~\cite{Ste}:

\begin{proposition} \cite{RR}\label{stein}
Let $(W,W')$ be an exchangeable pair so that $\mathbb{E}(W)=0$ and $\var(W)=1$ and assume
$$
\mathbb{E}(W'|W)=(1-\lambda)W
$$
for some $\lambda\in(0,1)$. Then for all real $t$:
$$
\left|\mathbb{P}(W\le t)-N(t)\right|\le\frac6\lambda\sqrt{\var\left(\mathbb{E}((W'-W)^2|W)\right)}
+6\sqrt{\frac1\lambda \mathbb{E}(|W'-W|^3)}.
$$
\end{proposition}

\noindent We proceed in five steps: (A) We begin with a 
classical `big block-small block' argument and approximate
$W_n=\frac{J_n}{\sigma_n}$ by a sum of random variables which
are separated by gaps. In (B) we then replace those random variables
by independent random variables. In (C) we define the interchangeable
pair in the usual way and estimate the terms on the right hand side of 
Proposition~\ref{stein}. In (D) and (E) we estimate the effects the steps (A) and
(B) have on the distributions.

We approximate $W_n=\frac{J_n}{\sigma_n}$ 
(clearly $\mathbb{E}(W_n)=0, \sigma(W_n)=1$)
by the random variable $\hat{W}_n=\frac1{\sqrt{r}}\sum_{j=0}^{r-1}W_m\circ T^{m'j}$ (that is
$\hat{W}_n=\frac1{\sqrt{r}\sigma_m}\sum_{j=0}^{r-1}J_m\circ T^{m'j}$)
where $m'=m+\Delta$ and $n=rm+(r-1)\Delta$. (For other values of $n$ not of this form we get 
an additional error term of the order $m'$.) \\
{\bf (A)} If we put $\hat{\cal A}^n=\bigvee_{j=0}^{r-1}T^{-m'j}{\cal A}^m$
then
\begin{eqnarray*}
\|\hat{W}_n-W_n\|_2&\le&\frac1{\sigma_n}\|J_{{\cal A}^n}-J_{\hat{\cal A}^n}\|_2
+\frac1{\sigma_n}\left\|I_{\hat{\cal A}^n}-\sum_{j=0}^{r-1}I_m\circ T^{m'j}\right\|_2\\
&&+\frac1{\sigma_n}|H(\hat{\cal A}^n)-rH_m|
+\left|\frac1{\sigma_n}-\frac1{\sqrt{r}\sigma_m}\right|\cdot\left\|\sum_{j=0}^{r-1}J_m\circ T^{m'j}\right\|_2.
\end{eqnarray*}
We individually estimate the four terms on the right hand side as follows: \\
{\bf (i)} By Lemma~\ref{variance.join} and Proposition~\ref{variance.convergence}
$$
\|J_{{\cal A}^n}-J_{\hat{\cal A}^n}\|_2=\sigma({\cal A}^n|\hat{\cal A}^n)
=\sigma\left(\left.\bigvee_{j=1}^{r-1}T^{-m-m'j}{\cal A}^\Delta\right|\hat{\cal A}^n\right)
=\sigma\left(\bigvee_{j=1}^{r-1}T^{-m-m'j}{\cal A}^\Delta\right)
\le c_1r\sqrt\Delta.
$$
{\bf (ii)} 
If  ${\cal D}_k=\bigvee_{j=0}^{k-1}T^{-m'j}{\cal A}^m$ then 
${\cal D}_{k+1}={\cal D}_k\vee T^{-m'k}{\cal A}^m$, $k=1,2,\dots,r$,
and by Lemma~\ref{rho.remainder.estimate} ($a=2$)
\begin{eqnarray*}
\|I_{{\cal D}_{k+1}}-I_{{\cal D}_k}-I_m\circ T^{m'k}\|_2^2
&=&\sum_{D\in{\cal D}_k,A\in T^{-m'k}{\cal A}^m}\mu(D\cap A)
\left(\frac1{\mu(D\cap A)}-\frac1{\mu(D)}-\frac1{\mu(A)}\right)^2\\
&=&\sum_{D\in{\cal D}_k,A\in T^{-m'k}{\cal A}^m}\mu(D\cap A)
\log^2\left(1+\frac{\rho(D,A)}{\mu(D)\mu(A)}\right)\\
&\le&c_2\left(\psi(\Delta)n^{3\beta}+n^{2\beta-(\beta-1)w}\right)
\end{eqnarray*}
for $k=1,2,\dots,r$. Hence (as ${\cal D}_1={\cal A}^m$)
$$
\|I_{\hat{\cal A}^n}-\sum_{j=0}^{r-1}I_m\circ T^{m'j}\|_2
\le\sum_{k=1}^r\|I_{{\cal D}_{k+1}}-I_{{\cal D}_k}-I_m\circ T^{m'k}\|_2
\le c_3r\sqrt{\psi(\Delta)n^{3\beta}+n^{2\beta-(\beta-1)w}}.
$$
{\bf (iii)} $|H(\hat{\cal A}^n)-rH_m|\le c_4r\left(\psi(\Delta)n^{2\beta}+n^{\beta-(\beta-1)w}\right)$
by Lemma~\ref{entropy.additivity}.\\
{\bf (iv)} Since by Proposition~\ref{variance.convergence} 
$$
\left|\frac1{\sigma_n}-\frac1{\sqrt{r}\sigma_m}\right|=\frac{|\sigma_n-\sqrt{r}\sigma_m|}{\sqrt{r}\sigma_n\sigma_m}
\le c_5\frac{m^{-\eta}}{\sqrt{n}},
$$
Lemma~\ref{variance.join} and again Proposition~\ref{variance.convergence} 
$$
\left\|\sum_{j=0}^{r-1}J_m\circ T^{m'j}\right\|_2
=\sigma\left(\bigvee_{j=0}^{r-1}T^{-m'j}{\cal A}^m\right)
\le r\sigma\left({\cal A}^m\right)
={\cal O}\left(r\sqrt{m}\right),
$$
we obtain that the fourth term is ${\cal O}(\sqrt{r}m^{-\eta})$, for any 
$\eta<\eta_0$.

\vspace{2mm}

\noindent Therefore, if $n$ is large enough,
\begin{eqnarray*}
\|\hat{W}_n-W_n\|_2&\le&c_6\left(\frac{r\sqrt\Delta}{\sqrt{n}}
+\frac{r}{\sqrt{n}}\sqrt{\psi(\Delta)n^{3\beta}+n^{2\beta-(\beta-1)w}}
+\frac{r}{\sqrt{n}}\left(\psi(\Delta)n^{2\beta}+n^{\beta-(\beta-1)w}\right)
+ \frac{\sqrt{r}}{m^{\eta}}\right)\\
&\le&c_7\left(\frac{r\Delta}{\sqrt{n}}
+rn^{\frac32\beta-\frac12}\sqrt{\psi(\Delta)}\left(1+n^{\frac12\beta}\sqrt{\psi(\Delta)}\right)
+rn^{\beta-\frac12-\frac12(\beta-1)w}+\frac{\sqrt{r}}{m^\eta}\right)
\end{eqnarray*}
as $\sigma_n\sim\sqrt{n}$ and $\beta-1>0$.

\vspace{3mm}

\noindent {\bf (B)}
Now let $X_j$ for $j=0,1,\dots,r-1$ be independent random variables
that have the same distributions as $W_m\circ T^{m'j}$, $j=0,1,\dots,r-1$. Put  $D_{V_n}(t)$
for the distribution function of the random variable $V_n=\frac1{\sqrt{r}}\sum_{j=0}^{r-1}X_j$ 
and $D_{\hat{W}_n}(t)$ for the distribution function of $\hat{W}_n$.
Since $V_n$ and $\hat{W}_n$ assume the same values, the difference between
the distributions is given by (with ${\cal D}_k=\bigvee_{j=0}^{k-1}T^{-m'j}{\cal A}^m$
as above):
\begin{eqnarray*}
\sup_t\left|D_{\hat{W}_n}(t)-D_{V_n}(t)\right|
&\le&\sum_{A_0\in {\cal A}^m}\cdots \sum_{A_{r-1}\in T^{-m'(r-1)}{\cal A}^m}
\left|\mu\left(\bigcap_jA_j\right)-\prod_j\mu(A_j)\right|\\
&\le&\sum_{k=0}^{r-1}\sum_{D\in{\cal D}_k}\sum_{A\in T^{-m'k}{\cal A}^m}
\left|\mu(D\cap A)-\mu(D)\mu(A)\right|\\
&=&\sum_{k=0}^{r-1}\sum_{D\in{\cal D}_k}\sum_{A\in T^{-m'k}{\cal A}^m}|\rho(D,A)|\\
&\le&c_8r\psi(\Delta)
\end{eqnarray*}
by the mixing property if we assume $n$ is large enough.

\vspace{3mm} 

\noindent {\bf (C)} In order to apply Proposition~\ref{stein} let us now define an
interchangeable pair in the  usual way by setting
$V'_n=V_n-\frac1{\sqrt{r}}X_Y+\frac1{\sqrt{r}}X^*$
where $Y\in\{0,1,\dots,r-1\}$ is a randomly chosen index and $X^*$ is a random
variable which is independent of all other random variables and has the same
distribution as the $X_j$. Since the random variables $X_j$ for $j=0,1,\dots,r-1$,
are i.i.d., the pair $(V'_n,V_n)$ is exchangeable. Moreover
$$
\mathbb{E}(V'_n|V_n)=\left(1-\frac1r\right)V_n
$$
(i.e. $\lambda=\frac1r$). 

We now estimate the two terms on the rights hand side of Proposition~\ref{stein}
separately:\\
{\bf (i)}
The third moment term of Proposition~\ref{stein} is estimated
using Corollary~\ref{third.moment}:
$$
\mathbb{E}(|V'_n-V_n|^3)^\frac13=\frac1{\sqrt{r}\sigma_m}\mathbb{E}(|J_m\circ T^{m'Y}+J_m^*|^3)^\frac13
\le\frac2{\sqrt{r}\sigma_m}M_3(J_m^3)^\frac13\le \frac{c_9}{\sqrt{r}}.
$$
Hence
$$
\sqrt{\frac1\lambda\mathbb{E}(|V'_n-V_n|^3)}
=\sqrt{\frac{c_9^3r}{r^\frac32}}
={\cal O}\left(r^{-\frac14}\right).
$$
{\bf (ii)} To estimate the variance term we follow Stein~\cite{Ste} and obtain
$$
\var\left(\mathbb{E}\left(\left.(V'_n-V_n)^2\right|V_n\right)\right)
\le\frac1{r^2}\var\left((X_Y-X^*)^2|X_0,X_1,\dots,X_{r-1}\right).
$$
Since
$$
\mathbb{E}(X_j^2|V_n)=\frac1r\sum_i\mathbb{E}(X_i^2|V_n)=\frac1r\sum_iX_i^2,
$$
we get
$$
\var(\mathbb{E}(X_Y^2|V_n))=\var\left(\frac1r\sum_iX_i^2\right)
=\frac1{r^2}\var\left(\sum_iX_i^2\right)=\frac1{r^2}r\var(X_0^2)=\frac1r\var(X_0^2).
$$
Since $X_0$ has the same distribution as $\frac1{\sigma_m}J_m$ we have 
$\mathbb{E}(X_0)=0$ and by Propositions~\ref{variance.convergence} and~\ref{fourth.moment}
$$
\var(X_0^2)=\var\left(\frac1{\sigma_m^2}J_m^2\right)=\frac1{\sigma_m^4}\sigma^2(J_m^2)
\le\frac1{\sigma_m^4}M_4({\cal A}^m)\le c_{10}.
$$
Hence 
$$
\frac6\lambda\sqrt{\var\left(\mathbb{E}((V'_n-V_n)^2|V_n)\right)}
\le c_{11}r\sqrt{\frac1{r^3}}\le c_{11}\frac1{\sqrt{r}}.
$$

Combining the estimates~(i) and~(ii) yields by Proposition~\ref{stein}
$$
\left|\mathbb{P}(V_n\le t)-N(t)\right|\le c_{11}\frac1{\sqrt{r}}+\frac{6\sqrt{c_9}}{r^\frac14}
\le c_{12}\frac1{\sqrt[4]{r}}.
$$
\vspace{3mm}

\noindent {\bf (D)} Part (B) and (C) combined yield
$$
\left|\mathbb{P}(\hat{W}_n\le t)-N(t)\right|
\le\left|\mathbb{P}(V_n\le t)-N(t)\right|+\left\|D_{\hat{W}_n}-D_{V_n}\right\|_\infty
\le  c_{12}\frac1{\sqrt[4]{r}}+c_8r\psi(\Delta).
$$
Let us put $\epsilon=\|W_n-\hat{W}_n\|_2$ and 
$\epsilon'=\sup_t\left|\mathbb{P}(\hat{W}_n\le t)-N(t)\right|$. Then 
($D_{W_n}$  is the distribution function of $W_n$) $N(t)\le \epsilon'$ for 
$t\le-|\log\epsilon'|$ and therefore $D_{\hat{W}_n}(t)\le2\epsilon'$ for $t\le-|\log\epsilon'|$
and similarly  \ $|1-N(t)|\le \epsilon'$ and consequently $|1-D_{\hat{W}_n}(t)|\le2\epsilon'$ 
for all $t\ge|\log\epsilon'|$ we get 
$$
\left\|(D_{W_n}-D_{\hat{W}_n})\chi_{[-|\log\epsilon'|,|\log\epsilon'|]}\right\|_\infty
\le2|\log\epsilon'|\cdot\left\|W_n-\hat{W}_n\right\|_2=2|\log\epsilon'|\epsilon
$$
and (since distribution functions are increasing)
$$
\left\|D_{W_n}-D_{\hat{W}_n}\right\|_\infty\le2|\log\epsilon'|\epsilon+2\epsilon'
$$

\vspace{3mm}

\noindent {\bf (E)} To optimise the bound
$$
\left|\mathbb{P}(W_n\le t)-N(t)\right|
\le \left|\mathbb{P}(\hat{W}_n\le t)-N(t)\right|+\|D_{W_n}-D_{\hat{W}_n}\|_\infty
\le 2|\log\epsilon'|\epsilon+3\epsilon'
$$
we distinguish between the case when $\psi$ decays (i) polynomially
and (ii) hyper polynomially.\\
{\bf (i)} Assume that $\psi$ decays polynomially with power $p>12$.
Let $\delta,\alpha\in(0,1)$ and put $m=[n^\alpha]$, $\Delta=[m^\delta]$ 
(i.e.\ $\Delta\sim n^{\alpha\delta}$, $\psi(\Delta)={\cal O}(n^{-\alpha\delta p}$).
Then (assuming $n^{\frac12\beta}\sqrt{\psi(\Delta)}={\cal  O}(1)$ 
 which will be satisfied once we choose $\beta$ and $\delta$)
$$
\|\hat{W}_n-W_n\|_2\le c_{13}\left(n^{\frac12-\alpha+\alpha\delta}
+n^{\frac12-\alpha+\frac32\beta-\frac12\alpha\delta p}
+n^{\frac12+\beta-\alpha-\frac12(\beta-1)w}+n^{\frac12-\frac\alpha2-\alpha\eta}\right).
$$
The first three terms on the right hand side are optimised by $\beta=\frac{w(p+2)}{(p+2)(w-2)+6}$
and $\alpha\delta=\frac{3\beta}{p+2}$. Then $\|\hat{W}_n-W_n\|_2\le\epsilon$,
$\epsilon= {\cal O}(n^x)$, where
$x=\max\left(\frac12-\alpha+\frac{3w}{(p+2)(w-2)+6},\frac12-\frac\alpha2-\alpha\eta\right)$.
The fourth term is smaller than the first three since we can assume that $\eta>\frac13$ as $w>4$.
The value of $\alpha$ is found by minimising the error term 
$2\epsilon|\log\epsilon'|+3\epsilon'$. Ignoring the logarithmic term we obtain
 $\alpha=\frac35+\frac{12}5\frac{w}{(p+2)(w-2)+6}$
which implies 
$$
\left|\mathbb{P}(W_n\le t)-N(t)\right|\le c_{14}\frac1{n^\kappa},
$$
for any $\kappa<\frac1{10}-\frac35\frac{w}{(p+2)(w-2)+6}$. Note
that $\alpha\eta>\kappa$ for all (possible) values of $p$ and $w$. \\
{\bf (ii)} If $\psi$ decays faster than 
any power then we can choose $\delta>0$ arbitrarily close to zero and
obtain $\alpha<\frac35$ which yields the estimate 
$
\left|\mathbb{P}(W_n\le t)-N(t)\right|\le c_{15}\frac1{n^\kappa},
$
for any $\kappa<\frac1{10}$. 

This concludes the proof since 
$W_n=\frac{I_n-H_n}{\sigma_n}$.
\qed

\vspace{3mm}

\noindent {\bf Proof of Theorem~\ref{CLT}.} We use Theorem~\ref{CLTn}  and have to make 
the following adjustments:\\
{\bf (i)} To adjust for the difference between $H_n$ and $nh$ we use
Lemma~\ref{entropy.approximation}:
$$
\mathbb{P}\left(\frac{I_n(x)-nh}{\sigma\sqrt{n}}\le t\right)
=\mathbb{P}\left(
\frac{I_n(x)-H_n}{\sigma\sqrt{n}}\le t+{\cal O}\left(n^{\frac12-\gamma}\right)\right)
=N(t)+{\cal O}\left(n^{-\kappa}\right)+{\cal O}\left(n^{\frac12-\gamma}\right).
$$
Since $p$ is big enough $\gamma$ can be chosen so that $\gamma-\frac12>\kappa$.\\
{\bf (ii)} By Proposition~\ref{variance.convergence} 
$\frac{\sigma_n}{\sqrt{n}}=\sigma+{\cal O}\left(n^{-\eta}\right)$ which yields
$$
\mathbb{P}\left(\frac{I_n(x)-H_n}{\sigma\sqrt{n}}\le t\right)
=\mathbb{P}\left(\frac{I_n(x)-H_n}{\sigma_n}\le t_n\right)
=N(t_n)+{\cal O}\left(n^{-\kappa}\right)
=N(t)+{\cal O}\left(n^{-\min(\eta,\kappa)}\right),
$$
where $t_n=t\frac{\sigma\sqrt{n}}{\sigma_n}=t\left(1+{\cal O}\left(n^{-\eta}\right)\right)$.
This concludes the proof since $\eta$ can be taken to be $>\kappa$.
 \qed

\section{Proof of Theorem~\ref{wip}  (Weak Invariance Principle)}
In order to prove the WIP for $I_n(x)=-\log\mu(A_n(x))$
denote by $W_{n,x}(t)$, $t\in [0,1]$, its interpolation
$$
W_{n,x}(k/n)=\frac{I_k(x)-kh}{\sigma\sqrt{n}}
$$
$x\in\Omega$ and linearly interpolated on
each of the subintervals $\left[\frac{k}n,\frac{k+1}n\right]$.
In particular $W_{n,x}\in C_\infty([0,1])$ (with supremum norm). Denote by $D_n$ the
distribution of $W_{n,x}$ on $C_\infty([0,1])$, namely
$$
D_n(H)=\mu\left(\left\{x\in\Omega: W_{n,x}\in H\right\}\right)
$$
where $H$ is a Borel subset of $C_\infty([0,1])$. The WIP then asserts
that the distribution $D_n$ converges weakly to the Wiener
measure, which means that $S_n=I_n-nh$ is for large $n$, and
after a suitable normalization distributed approximately as the
position at time $t=1$ of a particle in Brownian motion
\cite{Bi2}. 

If we put $S_i=-\log\mu(A_i(x))-ih(\mu)$
then two conditions have to be verified (\cite{Bi2} Theorem 8.1), namely\\
(A) The tightness condition: There exists a $\lambda>0$ so that for every
$\varepsilon>0$ there exists an $N_0$ so that
\begin{equation}\label{tightness}
\mathbb{P}\left(\max_{0\le i\le n}|S_i|>2\lambda\sqrt{n}\right)
\le\frac{\varepsilon}{\lambda^2}
\end{equation}
for all $n\ge N_0$.\\
\noindent (B) The finite-dimensional distributions of $S_i$ converge to those of the Wiener measure.

\vspace{2mm}

\noindent {\bf (A)} {\it Proof of tightness}:  As before let $J_i=I_i-H_i$ and note that
$ih-H_i={\cal O}(i^{1-\gamma})$, $1-\gamma\in(\frac{2w}{p(w-1)},1)$, 
(Lemma~\ref{entropy.approximation})
 is easily absorbed by the term $\lambda\sqrt{n}$ as $1-\gamma<\frac12$.
In the usual way (cf.\ e.g.\
\cite{Bi2}) we get
$$
\mathbb{P}\left(\max_{0\le i\le n}| J_i|>2\lambda\sqrt{n}\right)
\le\mathbb{P}\left(|J_n|>\lambda\sqrt{n}\right)
+\sum_{i=0}^{n-1}\mu\left(E_i\cap\{|J_i-J_n|\ge\lambda\sqrt{n}\}\right),
$$
where $E_i$ is the set of points $x$ so that
$|J_i(x)|>2\lambda\sqrt{n}$ and $|J_k(x)|\le2\lambda\sqrt{n}$ for
$k=0,\dots,i-1$. Note that $E_i$ lies in the $\sigma$-algebra
generated by ${\cal A}^i$. Clearly the sets $E_i$ are pairwise
disjoint. 
To estimate $\mu\left(E_i\cap\{|J_i-J_n|\ge\lambda\sqrt{n}\}\right)$
let us first `open a gap' of length $\Delta<\frac{n}2$. Let 
$\tilde{\cal A}^n={\cal A}^i\vee T^{-i-\Delta}{\cal A}^{n-i-\Delta}$
(if $i<\frac{n}2$ and $\tilde{\cal A}^n={\cal A}^{i-\Delta}\vee T^{-i}{\cal A}^{n-\Delta}$
if $i\ge\frac{n}2$),
denote by $\tilde{I}_n$ its information function and by $\tilde{H}_n=\mu(\tilde{I}_n)$
its entropy. Obviously $H_n\ge\tilde{H}_n$ and moreover
$\mu(I_n-\tilde{I}_n)=H_n-\tilde{H}_n\le H_\Delta\le c_1\Delta$. 
Since by Lemma~\ref{variance.join} and Corollary~\ref{K.estimate} (as ${\cal A}^n$
refines $\tilde{\cal A}^n$) 
$$
\sigma(I_n-\tilde{I}_n)=\sigma({\cal A}^n|\tilde{\cal A}^n)
\le\sqrt{K_2({\cal A}^\Delta)}\le c_2\Delta
$$
we obtain by Chebycheff's inequality ($\tilde{J}_n=\tilde{I}_n-\tilde{H}_n$)
\begin{equation}\label{first}
\mathbb{P}(|J_n-\tilde{J}_n|\ge\ell)\le\frac{\sigma^2(I_n-\tilde{I}_n)}{\ell^2}\le c_3\frac{\Delta^2}{\ell^2}.
\end{equation}
By the uniform strong mixing property
$$
\tilde{I}_n(B)=I_i(B)+I_{n-i-\Delta}(C)-\log\left(1+\frac{\rho(B,C)}{\mu(B)\mu(C)}\right)
$$
for all $(B,C)\in{\cal A}^i\times T^{-i-\Delta}{\cal A}^{n-i-\Delta}$.
If $Y$ denotes the random variable on ${\cal A}^i\times T^{-i-\Delta}{\cal A}^{n-i-\Delta}$
whose values are $Y(B,C)=-\log\left(1+\frac{\rho(B,C)}{\mu(B)\mu(C)}\right)$ then
by Lemma~\ref{rho.remainder.estimate} ($a=2$)
$$
\sigma^2(Y)\le\|Y\|_{L^2}^2
\le C_4\left(\psi(\Delta)(n-\Delta)^{3\beta}+(n-\Delta)^{2\beta-(\beta-1)w}\right)
$$
for $\beta>1$ arbitrary. By Chebycheff's inequality this implies
\begin{equation}\label{second}
\mathbb{P}\left(\left|\tilde{J}_n-J_i-J_{n-i-\Delta}\circ T^{i+\Delta}\right|\ge\ell\right)
\le\frac{\sigma^2(Y)}{\ell^2}
\le C_4\frac{\psi(\Delta)(n-\Delta)^{3\beta}+n^{2\beta-(\beta-1)w}}{\ell^2}.
\end{equation}
Then
\begin{eqnarray*}
\mu\left(E_i\cap\left\{|J_n-J_i|\ge\lambda\sqrt{n}\right\}\right)
&\le&\mu\left(E_i\cap\left\{|J_n-\tilde{J}_n|\ge\ell\right\}\right)
+\mu\left(E_i\cap\left\{|\tilde{J}_n-J_i-J_{n-i-\Delta}\circ T^{i+\Delta}|\ge\ell\right\}\right)\\
&&\hspace{3cm}+\mu\left(E_i\cap\left\{|J_{n-i-\Delta}\circ T^{i+\Delta}|\ge\lambda\sqrt{n}-2\ell\right\}\right).
\end{eqnarray*}
The last term on the right hand side can be estimated using the mixing 
property (note that $E_i$ is in the $\sigma$-algebra generated by ${\cal A}^i$,
and $\{|J_{n-i-\Delta}|\ge\lambda\sqrt{n}-2\ell\}$ is in the $\sigma$-algebra generated by
$T^{-i-\Delta}{\cal A}^{n-i-\Delta}$)
\begin{eqnarray*}
\mu\left(E_i\cap\left\{|J_{n-i-\Delta}\circ T^{i+\Delta}|\ge\lambda\sqrt{n}-2\ell\right\}\right)
&=&\mu(E_i)\mathbb{P}\left(|J_{n-i-\Delta}|\ge \lambda\sqrt{n}-2\ell\right)\\
&&+\sum_{B\subset E_i}\;\;\sum_{C\subset T^{-i-\Delta}\{|J_{n-i-\Delta}|\ge\lambda\sqrt{n}-2\ell\}}\rho(B,C)\\
&\le&\mu(E_i)\left(2N\left(\frac{\lambda\sqrt{n}-2\ell}{\sigma_{n-i-\Delta}}\right)
+C_0(n-i-\Delta)^{-\kappa}\right)+\psi(\Delta)
\end{eqnarray*}
using Theorem~\ref{CLTn} in the last step.

We finally obtain
(as $\mathbb{P}(|J_n|>\lambda\sqrt{n})\le2N(\lambda)+c_4n^{-\kappa}$)
\begin{eqnarray*}
\mathbb{P}\left(\max_{0\le i\le n}|J_i|>2\lambda\sqrt{n}\right)
&\le&2N(\lambda)+c_4n^{-\kappa}
+\sum_i\mu\left(E_i\cap\{|J_n-\tilde{J}_n|\ge\ell\}\right)\\
&&+nC_4\frac{\psi(\Delta)n^{3\beta}+n^{2\beta-(\beta-1)w}}{\ell^2}\\
&&+\sum_i\mu(E_i)\left(2N\left(\frac{\lambda\sqrt{n}-2\ell}{\sigma_{n-i-\Delta}}\right)
+C_0(n-i-\Delta)^{-\kappa}\right)+n\psi(\Delta)\\
&\le&2N(\lambda)+c_5n^{-\kappa}
+c_6\frac{\Delta^2+\psi(\Delta)n^{3\beta}+n^{2\beta-(\beta-1)w}}{\ell^2}
+2N\left(\frac{\lambda\sqrt{n}-2\ell}{\sqrt{n}}\right)
\end{eqnarray*}
(if $\Delta<\frac{n}2$ is small enough). If $\psi$ decays at least polynomially with a power
larger than $8+\frac{24}{w-4}$ then we can put $\ell\sim n^\alpha, \Delta\sim n^{\alpha'}$
and choose $\alpha'<\alpha<\frac12$ and $\beta>1$ (e.g.\ $\beta=\frac{w}{w-2}$,
$\alpha'<\frac{3\beta}p$) so that the terms on the right hand side
which don't involve the normal probability $N$ decay polynomially in $n$.
This proves the tightness condition (\ref{tightness}), since
for every $\varepsilon>0$ one can find a $\lambda>1$ so that
the quadratic estimate holds for all $n$ large enough. 

\vspace{2mm}

\noindent {\bf (B)} {\it Proof of the finite-dimensional distribution
convergence:} For $t\in[0,1]$ define the random variable
$$
X_n(t,x)=\frac1{\sigma\sqrt{n}}\left(S_{[nt]}(x)
+(nt-[nt])\left(S_{[nt]+1}(x)-S_{[nt]}(x)\right)\right)
$$
which interpolates $S_{[nt]}$. It is defined on $\Omega$ and has values in 
$C_\infty([0,1])$.

We must show that the distribution of
$(X_n(t,x),X_n(t,x)-X_n(s,x))$ converges
to $({\cal N}(0,t),{\cal N}(0,t-s))$ ($0\le s<t$) as $n\rightarrow\infty$,
where ${\cal N}(0,t)$ is the normal distribution with zero mean
and variance $t^2$.
To prove this as well as the convergence of higher finite
dimensional distributions it suffices to show that $X_n(t,x)-X_n(s,x)$ 
 converges to ${\cal N}(0,t-s)$  (\cite{Bi2} Theorem 3.2).
We obtain by  Lemma~\ref{entropy.approximation}
\begin{eqnarray*}
S_{[nt]}-S_{[ns]}&=&J_{[nt]}-J_{[ns]}+{\cal O}\left((nt\right)^{1-\gamma})
\end{eqnarray*}
and by (\ref{first}), (\ref{second})  and Theorem~\ref{CLT}
\begin{eqnarray*}
\mathbb{P}\left(\frac{S_{[nt]}-S_{[ns]}}{\sigma\sqrt{n}}\ge\lambda\right)
&\le&\mathbb{P}\left(\left|J_{[nt]}-\tilde{J}_{[nt]}\right|\ge\ell\right)
+\mathbb{P}\left(\left|\tilde{J}_{[nt]}-J_{[ns]}-J_{[nt]-[ns]-\Delta}\circ T^{[ns]+\Delta}\right|\ge\ell \right)\\
&&+\mathbb{P}\left(\left|J_{[nt]-[ns]-\Delta}\right|\ge\lambda\sigma\sqrt{n}-2\ell \right)
+{\cal O}\left((nt)^{\frac12-\gamma}\right)\\
&\le&\frac{\sigma^2(I_n-\tilde{I}_n)}{\ell^2}+\frac{\sigma^2(Y)}{\ell^2}
+N\left(\frac{\lambda\sigma\sqrt{n}-2\ell}{\sqrt{[nt]-[ns]-\Delta}}\right)
+\frac{C_0}{([nt]-[ns]-\Delta)^{\kappa}}+\frac{{\cal O}(1)}{(nt)^{1-\gamma}}\\
&\le&c_3\frac{\Delta^2}{\ell^2}+C_4\frac{\psi(\Delta)(nt)^{3\beta}+(nt)^{2\beta-(\beta-1)w}}{\ell^2}
+\frac{c_7}{(n(t-s))^{\kappa}}+N\left(\frac{\lambda}{\sqrt{t-s}}\right),
\end{eqnarray*}
assuming $\frac12-\gamma\ge\kappa$ and $n(t-s)>\!\!>\Delta$. Similarly to above we used
a random variable $Y$ on ${\cal A}^{[ns]}\times T^{-[ns]-\Delta}{\cal A}^{[nt]-[ns]-\Delta}$
given by $Y(B,C)=-\log\left(1+\frac{\rho(B,C)}{\mu(B)\mu(C)}\right)$. 
Now let $\ell\sim n^\alpha, \Delta\sim n^{\alpha'}$
and $\alpha'<\alpha<\frac12$ and $\beta>1$ so that the terms on the right hand side
other than $N\left(\lambda/\sqrt{t-s}\right)$ decay polynomially in $n$.
Hence $S_{[nt]}-S_{[ns]}$ and therefore $X_n(t,x)-X_n(s,x)$ converges in distribution to
${\cal N}(0,\sqrt{t-s})$ as $n\rightarrow\infty$. \qed

\section{Appendix (Markov chains)}

Here we compute the variance for the Markov measure on an infinite alphabet.
As in section~\ref{examples} let $\Sigma$ be the shiftspace over the alphabet
$\mathbb{N}$ and $\mu$ the Markov measure generated by the probability vector
$\vec{p}$ and stochastic matrix $P$. 
Then
$$
\sigma_n^2=\frac12\sum_{\vec{x},\vec{y}\in{\cal A}^n}\mu(\vec{x})\mu(\vec{y})
\left(\log\frac{p_{x_1}}{p_{y_1}}+\sum_{j=1}^{n-1}\log\frac{P_{x_jx_{j+1}}}{P_{y_jy_{j+1}}}\right)^2
=A+B+C+D,
$$
where 
$$
A=\frac12\sum_{\vec{x},\vec{y}\in{\cal A}^n}\mu(\vec{x})\mu(\vec{y})\log^2\frac{p_{x_1}}{p_{y_1}}
=\frac12\sum_{ij}p_ip_j\log^2\frac{p_i}{p_j}={\cal O}(1)
$$
and 
\begin{eqnarray*}
B&=&\sum_{j=1}^{n-1}\sum_{\vec{x},\vec{y}\in{\cal A}^n}\mu(\vec{x})\mu(\vec{y})
\log\frac{p_{x_1}}{p_{y_1}}\log\frac{P_{x_jx_{j+1}}}{P_{y_jy_{j+1}}}\\
&=&\sum_{j=1}^{n-1}\sum_{\vec{x},\vec{y}\in{\cal A}^{j+1}}\mu(\vec{x})\mu(\vec{y})
\left(\log p_{x_1}\log P_{x_jx_{j+1}}+\log p_{y_1}\log P_{y_jy_{j+1}}
-\log p_{x_1}\log P_{y_jy_{j+1}}-\log p_{y_1}\log P_{x_jx_{j+1}}\right)\\
&=&2\sum_{j=1}^{n-1}\sum_{\vec{x}\in{\cal A}^{j+1}}\mu(\vec{x})\log p_{x_1}\log P_{x_jx_{j+1}}
+2(n-1)h\sum_ip_i\log p_i.
\end{eqnarray*}
Since Markov chains are exponentially mixing~\cite{Bre}
 we get for some $\vartheta\in(0,1)$ that 
$$
\sum_{\vec{x}\in{\cal A}^{j+1}}\mu(\vec{x})\log p_{x_1}\log P_{x_jx_{j+1}}
=\sum_ip_i\log p_i\sum_{ij}p_iP_{ij}\log P_{ij}+{\cal O}(\vartheta^j)
=-h\sum_ip_i\log p_i+{\cal O}(\vartheta^j)
$$
and therefore
$$
B=2\sum_j{\cal O}(\vartheta^j)={\cal O}(1).
$$
The principal term is
$$
D=\frac12\sum_{j=1}^{n-1}\sum_{\vec{x},\vec{y}\in{\cal A}^n}\mu(\vec{x})\mu(\vec{y})
\log^2\frac{P_{x_jx_{j+1}}}{P_{y_jy_{j+1}}}
=\frac{n-1}2\sum_{ijk\ell}p_iP_{ij}p_kP_{k\ell}\log^2\frac{P_{ij}}{P_{k\ell}}.
$$
Lastly we get the correction term
\begin{eqnarray*}
C&=&\sum_{i\not=j}\sum_{\vec{x},\vec{y}\in{\cal A}^n}\mu(\vec{x})\mu(\vec{y})
\log\frac{P_{x_ix_{i+1}}}{P_{y_iy_{i+1}}}\log\frac{P_{x_jx_{j+1}}}{P_{y_jy_{j+1}}}\\
&=&2\sum_{k=1}^{n-1}(n-k)\sum_{\vec{x},\vec{y}\in{\cal A}^{k+1}}\mu(\vec{x})\mu(\vec{y})
\log\frac{P_{x_1x_2}}{P_{y_1y_2}}\log\frac{P_{x_kx_{k+1}}}{P_{y_ky_{k+1}}}\\
&=&2\sum_{k=1}^{n-1}(n-k)\sum_{\vec{x},\vec{y}\in{\cal A}^{k+1}}\mu(\vec{x})\mu(\vec{y})
\left(\log P_{x_1x_2}\log P_{x_kx_{k+1}}+\log P_{y_1y_2}\log P_{y_ky_{k+1}}
-\log P_{x_1x_2}\log P_{y_ky_{k+1}}-\log P_{y_1y_2}\log P_{x_kx_{k+1}}\right)\\
&=&4\sum_{k=1}^{n-1}(n-k)\left(\sum_{\vec{x}\in{\cal A}^{k+1}}\mu(\vec{x})
\log P_{x_1x_2}\log P_{x_kx_{k+1}}-h^2\right).
\end{eqnarray*}
Since $\sigma^2=\lim_{n\rightarrow\infty}\frac{\sigma_n^2}n$
we finally obtain
$$
\sigma^2=\frac12\sum_{ijk\ell}p_iP_{ij}p_kP_{k\ell}\log^2\frac{P_{ij}}{P_{k\ell}}
+4\sum_{k=1}^\infty\sum_{\vec{x}\in{\cal A}^{k+1}}\mu(\vec{x})
\left(\log P_{x_1x_2}\log P_{x_kx_{k+1}}-h^2\right),
$$
where the infinite sum converges because the terms (correlations) decay exponentially fast.


\end{document}